\documentclass[onefignum,onetabnum,a4paper]{article}

\usepackage{graphicx}
\usepackage{amsmath,amsfonts,amssymb}  
\usepackage{url}  
\usepackage{psfrag}
\usepackage{sidecap}
\usepackage{colortbl}
\usepackage{color}
\usepackage{booktabs}
\usepackage{mathtools}
\usepackage{xcolor}
\usepackage{hyperref}
\usepackage[all]{hypcap}
\usepackage{multirow}
\usepackage{soul}
\usepackage[shortlabels]{enumitem}  
\usepackage{mathtools}

\definecolor{g}   {rgb}{0.7 0.7 0.6}  
\definecolor{r}   {rgb}{1 0.0 0.0}  

\begin{document}

\markboth{R. Otupiri, B. Krauskopf, N.G.R. Broderick}{The Yamada model for a self-pulsing laser with non-identical decay times}

\pagestyle{myheadings}

\title{The Yamada model for a self-pulsing laser: bifurcation structure for non-identical decay times of  gain and absorber}

\author{Robert Otupiri$^{\ddagger}$\footnote{{\tt r.otupiri@auckland.ac.nz}},
Bernd Krauskopf$^{\dagger}$\footnote{{\tt b.krauskopf@auckland.ac.nz}}, and 
Neil G.R. Broderick$^{\ddagger}$\footnote {n.broderick@auckland.ac.nz} \\
Department of Physics$^{\ddagger}$ and Department of Mathematics$^{\dagger}$,\\ The Dodd-Walls Centre for Photonic and Quantum Technologies,
\\The University of Auckland, Auckland 1142, New Zealand}

\date{\today}
\maketitle

\begin{abstract}
We consider self-pulsing in lasers with a gain section and an absorber section via a mechanism known as Q-switching, as described mathematically by the Yamada ordinary differential equation model for the gain, the absorber and the laser intensity. More specifically, we are interested in the case that gain and absorber decay on different time scales. We present the overall bifurcation structure by showing how the two-parameter bifurcation diagram in the plane of pump strength versus decay rate of the gain changes with the ratio between the two decay rates. In total, there are ten cases BI to BX of qualitatively different two-parameter bifurcation diagrams, which we present with an explanation of the transitions between them. Moroever, we show for each of the associated eleven cases of structurally stable phase portraits (in open regions of the parameter space) a three-dimensional representation of the organisation of phase space by the two-dimensional manifolds of saddle equilibria and saddle periodic orbits.

The overall bifurcation structure constitutes a comprehensive picture of the exact nature of the observable dynamics, including multi-stability and excitability, which we expect to be of relevance for experimental work on Q-switching lasers with different kinds of saturable absorbers.
\end{abstract}

\section{Introduction}

Pulsating oscillations in lasers have been observed in original studies in the 1960s; see, for example,~\cite{abraham1988dynamical}. Subsequently, pulsating output light was found to be present in a plethora of laser systems, including gas lasers~\cite{velikhov1983pulsed,heard1963ultra,little1999metal}, solid state lasers~\cite{magni1986resonators,yu2013topological} and semiconductor lasers~\cite{keller1996semiconductor,rulliere2005femtosecond}. Pulsing lasers have been used in diverse applications, such as telecommunications~\cite{tam1986applications}, millimetre wave generation to provide increased communication bandwidth~\cite{novak1994millimetre}, photo-acoustic sensing for increased accuracy~\cite{duarte1990dye,hamlin1991high} and medical imaging~\cite{gibson1993lasers}. A comprehensive review of pulsing lasers and their applications is provided, for example, by~\cite{svelto1998principles}.

A self-pulsing laser is an economical and simple means to provide short, high power laser pulses~\cite{siegman1971introduction}. In a much used approach, two spatially separated segments, a gain and a saturable absorber section, are placed in the laser cavity. These sections perform different roles in the system to obtain the desired output. The two materials are separated, so that they interact only through the intensity field in the cavity~\cite{erneux1988q}. The gain section is pumped so that the atoms have a positive population inversion consisting of excited states that provide an active or amplifying medium capable of producing light by producing photons. The saturable absorber, on the other hand, provides a negative population inversion initially, that is, it absorbs photons or light rather than producing them; this means that the absorber section appears opaque to the light inside the cavity rather than transparent. However, as time passes, the absorber gradually saturates, meaning that it absorbs less and less light, and eventually becomes fully saturated and transparent. This enables a short pulse to from and escape the laser cavity as a spike in the laser intensity. As a result, the absorber is flushed out and becomes opaque again. In the self-pulsing regime, this process, which is referred to as Q-switching, then repeats. The self-pulsing operation in such a laser with saturable absorber, hence, consists of a sequence of high-intensity pulses with the laser effectively being off in between the pulses. We remark that Q-switching requires the laser to operate in a single mode of the laser cavity. An alternative and physically different method for pulse generation is mode locking, which relies on the technique of locking the phases of many longitudinally propagating modes within the cavity~\cite{haus2000mode,smith1970mode}; as such, it is not the focus of this study.

Several models were developed for Q-switched single mode laser systems; they describe self-pulsations as well as interesting multi-stable, excitable and chaotic dynamics~\cite{yamada_theoretical_1993,ueno1985conditions,erneux1988q}. In this paper, we study a model for a self-pulsing laser, introduced by Yamada~\cite{yamada_theoretical_1993}, which comprises three ordinary differential equations (ODEs) for the gain G, the absorption Q and the intensity I. This model can be written in non-dimensional form~\cite{Dubbeldam_Krauskopf_Self_pulsations_lasers_saturable_absorber} as

\begin{align}
\centering
\label{Yamadamodel1}
\begin{split}
\dot{G} &= \gamma_{G} (A-G-GI),
\\
\dot{Q} &= \gamma_{Q} (B-Q-\frac{\gamma _{Q}}{\gamma _{G}} aQI),
\\
\dot{I} &= (G-Q-1)I.
\end{split}
\end{align}

\noindent Here, $A$ is the pump parameter for the gain, which is typically varied in an experiment. Further, $\gamma_{G}$ and $\gamma_{Q}$ are the ratios of  the photon versus carrier lifetimes, or decay times for short, of the gain and absorber sections, respectively, $B$ is the absorption coefficient, and $a$ the relative absorption versus gain parameter. The latter group of parameters is determined by the material properties of the laser device. We rewrite Eqs.~\eqref{Yamadamodel1} in terms of the decay time ratio $\sigma=\gamma _{G}/\gamma _{Q}$ as
\begin{align}
\centering
\label{Yamadamodel2}
\begin{split}
\dot{G} &=  \gamma_{G} (A-G-GI),
\\
\dot{Q} &= \frac{\gamma_{G}}{ \sigma} (B-Q-\frac{a}{\sigma}QI),
\\
\dot{I} &= (G-Q-1)I. 
\\
\end{split}
\end{align}
The new parameter $\sigma$ encodes the different decay times of gain and absorber media and, hence, their different decay rates. 
 
The Yamada model has been used by several authors to study Q-switching dynamics. In~\cite{peterson1999dynamics,erneux2000pulse}, it is suggested that Q-switched pulsing solutions are born from a homoclinic bifurcation; in addition an asymptotic method is used to map out the threshold conditions and pulse characteristics in a microchip laser with saturable absorber as a function of the pump parameter. The paper \cite{Dubbeldam_Krauskopf_Self_pulsations_lasers_saturable_absorber} presented a complete description of all possible dynamics for the special case that $\gamma_{G} = \gamma_{Q} = \gamma$, that is, for $\sigma=1$. More specifically, these authors presented a two-parameter bifurcation diagram in plane of $(A,\gamma)$-plane with a total of nine regions of different qualitative dynamics, including regions with bistability and excitability; moreover, they showed that self-pulsations are indeed born from a homoclinic bifurcation. That work was used in~\cite{barbay2011excitability,selmi2016spike} to explain the short pulses in experiments on micropillar lasers with an intracavity saturable absorber in response to an input perturbation. \cite{shastri2013graphene} demonstrated a unified platform for spike processing with a fibre laser with graphene saturable absorber; aided by the Yamada model, they showed that this platform can simultaneously exhibit logic-level restoration, cascadability and input-output isolation. More recently, in \cite{otupiri2018experimental} we considered a relatively short all-fibre laser with saturable absorber and demonstrated that the onset of self-pulsations is either via a homoclinic bifurcation or a Hopf bifurcation; these experimental results are in good agreement with the predictions of the Yamada model.

The motivation for this study stems from the fact that most experimental realisations of self-pulsing lasers have different decay times for gain and absorber sections and, hence, $\sigma \neq 1$. On the other hand, the earlier bifurcation analysis of the Yamada model in~\cite{Dubbeldam_Krauskopf_Self_pulsations_lasers_saturable_absorber} is restricted to the special case $\sigma = 1$. This discrepancy is the central issue that this paper addresses. More specifically, we present here a three-parameter bifurcation analysis by considering how the two-parameter bifurcation diagram in the $(A,\gamma_G)$-plane of Eqs.~\eqref{Yamadamodel2} changes with the decay ratio parameter $\sigma$. There is a total of ten cases, which we denote BI to BX, of qualitatively different two-parameter bifurcation diagrams, with a total of eleven cases of structurally stable phase portraits.
We present careful numerical evidence for all of the cases BI to BX and the transitions between them. Moreover, for each of the eleven phase portraits we present three-dimensional representation of the organisation of phase space, rather than two-dimensional projections as in~\cite{Dubbeldam_Krauskopf_Self_pulsations_lasers_saturable_absorber}. To this end, we compute and show the two-dimensional invariant manifolds of the respective saddle equilibria and saddle periodic orbits. 

Our first result is that the two-parameter bifurcation diagram from~\cite{Dubbeldam_Krauskopf_Self_pulsations_lasers_saturable_absorber} with a Bogdanov-Takens (BT) point as an organising center and featuring bistability and excitability --- which is our starting point and denoted case BI here --- does not change qualitatively when $\sigma$ is allowed to change for $\sigma < 1$. Hence, when the gain decays slower than the absorption in the respective sections of the laser, the $(A,\gamma_{G})$-plane is already known from this earlier work. On the other hand, as we will show, for $1 < \sigma$, when the gain decays faster than the absorption, there are qualitative changes to the $(A,\gamma_{G})$-plane from a certain value of $1 < \sigma$ onwards, which generate an additional two generic phase portraits. As our main result, we explain in detail how and when the bifurcation diagram in $(A,\gamma_{G})$-plane changes as $\sigma$ increases, by providing sketches of all bifurcation cases BI--BX as well as explanations of the transitions between them via codimension-three events. These sketches are accompanied by evidence in the $(A,\gamma_{G})$-plane of careful numerical investigations, including enlargements of regions that verify the presence or absence of certain dynamics. We show that the bifurcation diagram gradually decreases in size and eventually collapses down to a point with $\gamma_{G}=0$ as $\sigma$ increases. 

An overview of how the loci of the different bifurcations change with $\sigma$ is provided in two ways. Firstly, we present the $\sigma$-line with ranges of all bifurcation cases BI--BX and the locations of associated transitions, both as a sketch and as an image with the actual positions of numerically computed transitions. Secondly, we provide an image of the bifurcation set in three-dimensional $(A,\gamma_{G},\sigma)$-space, consisting of surfaces of codimension-one bifurcations, their intersection curves, as well as curves of codimension-two bifurcations. All eleven cases of generic phase portraits in open regions of $(A,\gamma_{G},\sigma)$-space are presented as two-dimensional sketches and projections, but also in the full $(G, Q, I)$-space. The latter provides new insight into how the two-dimensional invariant manifolds of equilibria and periodic orbits of saddle type form the basin boundaries between different attractors. In this way, the organisation of multi-stabilty, as well as the exact nature of the excitability threshold in the full three-dimensional phase space of this laser system, are clarified.

\section{Bifurcation diagram BI for $\sigma \leq 1$}
\label{sec:Bifurcation Analysis}

\begin{figure}[t!]
\centering
\includegraphics[width=0.8\columnwidth]{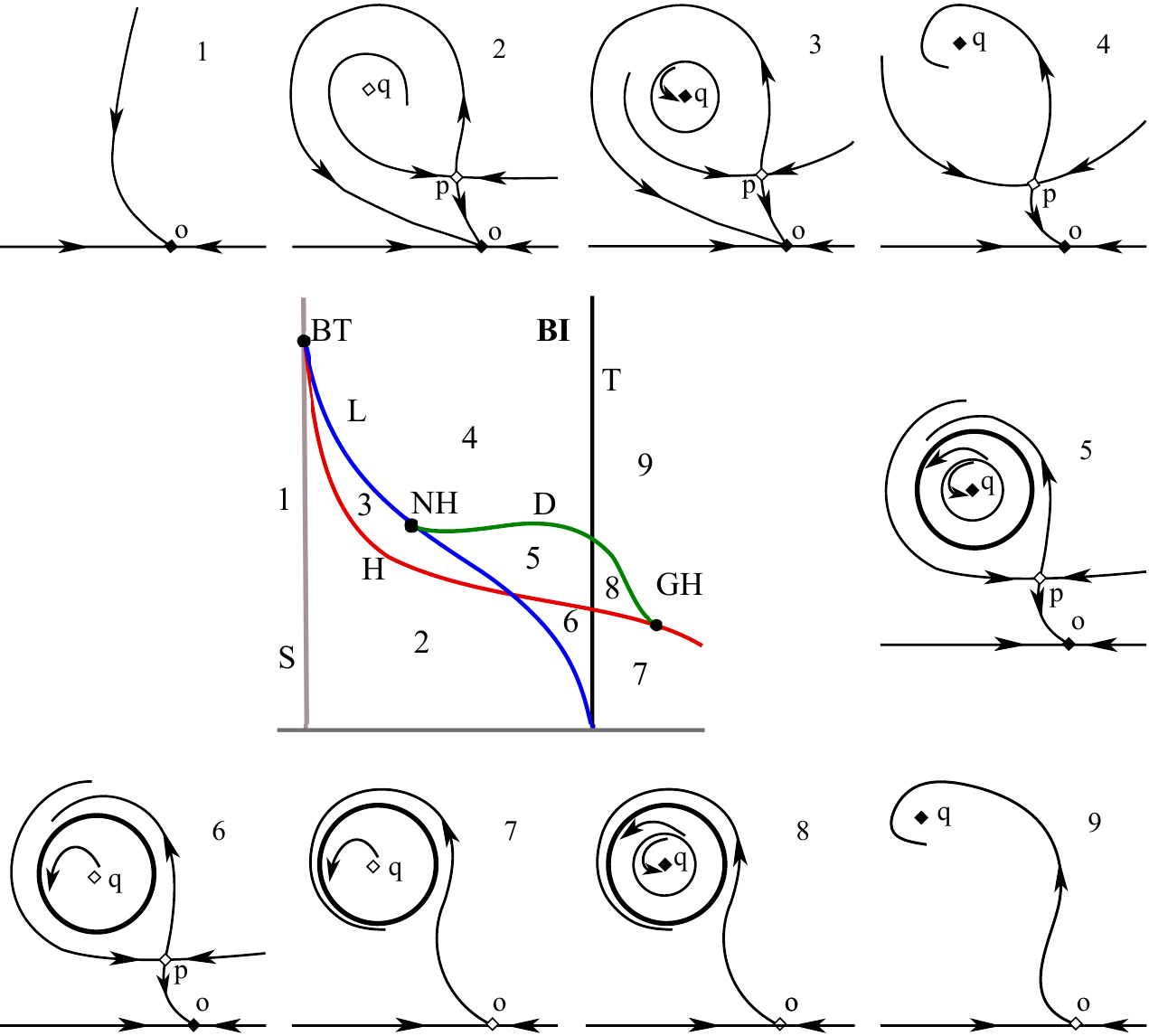}
\caption{Sketches of the bifurcation diagram of type \textbf{BI} in the $(A,\gamma)$-plane with phase portraits for the nine open regions, shown in the $(G, I)$-plane where the third direction is always attracting. Invariant objects are attracting equilibria (filled diamond), saddle equilibria (hollow diamond), attracting periodic orbits (boldface circle) and saddle periodic orbits (thin circles).}
\label{Bifurcation_Sketches_original_study}
\end{figure}

Our starting point is the detailed study in~\cite{Dubbeldam_Krauskopf_Self_pulsations_lasers_saturable_absorber} of self pulsations of a laser with saturable absorber described by the Yamada model for equal decay lifetimes, that is, for $\gamma=\gamma_{G} = \gamma_{Q}$ and, hence, $\sigma = 1$. A complete description of how the possible dynamics of the model depends on all four parameters $A$, $\gamma_{G}$, $B$  and $a$ was presented. Figure~\ref{Bifurcation_Sketches_original_study} shows the main bifurcation diagram from the paper~\cite{Dubbeldam_Krauskopf_Self_pulsations_lasers_saturable_absorber}, which includes an explicit derivation of all bifurcations of equilibria. A two-parameter bifurcation diagram in the $(A,\gamma_{G})$-plane is surrounded by sketches of nine phase portraits in open regions. This most complicated and relevant bifurcation diagram occurs when $B$ and $a$ are sufficiently large as is the case, for example, for $B=5.8$ and $a=1.8$, which are the values that we use for the computations throughout this paper. This bifurcation diagram is referred to in~\cite{Dubbeldam_Krauskopf_Self_pulsations_lasers_saturable_absorber} as Type III.

Since it is our starting point, we refer to this bifurcation diagram in Fig.~\ref{Bifurcation_Sketches_original_study} as case BI. The horizontal line corresponds to $\gamma_{G}=0$, and it is the lower limit of the bifurcation diagram because negative  $\gamma_{G}$ is not physically relevant. There are curves of saddle-node bifurcations S, and of transcritical bifurcation T. On the curve S there is a Bogdanov-Takens point BT, which emerges as the main organising centre; for background on the different types of bifurcations we encounter see, for example, \cite{GH,Kuz}. From BT, a curve H of Hopf bifurcations and a curve L of homoclinic bifurcations emerge. They cross one another in a codimension-one-plus-one event, H$\,\cap\,$L, and L stays on the left of the transcritical curve T where the equilibrium of the homoclinic orbit disappears. The curve H, on the other hand, crosses the curve T and eventually reaches $\gamma_{G}=0$ for higher values of  $A$ (which is not shown in the bifurcation sketch) of Fig.~\ref{Bifurcation_Sketches_original_study}. The crossing of curves L and H requires a switch of the stability of the bifurcating periodic orbits, that is, there must be changes of the criticality of L and H near their crossing point~\cite{Dubbeldam_Krauskopf_Self_pulsations_lasers_saturable_absorber}. This happens on the curve H at the codimension-two degenerate Hopf point GH where H changes from subcritical to supercritical; and on the curve L we find a codimension-two neutral-saddle homoclinic point NH, where the saddle quantity at the saddle point involved in the homoclinic bifurcation is zero; this means that the saddle is neither attracting nor repelling along the homoclinic orbit. Consequently, there is a curve D of saddle-node bifurcations of periodic orbits emanating from the points NH and GH, which connects these two points and completes the local picture near the intersection of L and H.

This collection of bifurcation curves divides the $(A,\gamma_{G})$-plane in  Fig.~\ref{Bifurcation_Sketches_original_study} into nine distinct regions of different dynamics, which are represented by the surrounding setches of the associated phase portraits 1 to 9. These are represented in two-dimensions only, where $I$ is the vertical direction with $I=0$ being the bottom line. The horizontal direction can be thought of either $G$ or $Q$ and the third direction is always attracting.

\begin{figure}[t!]
\centering
\includegraphics[width=0.95 \columnwidth]{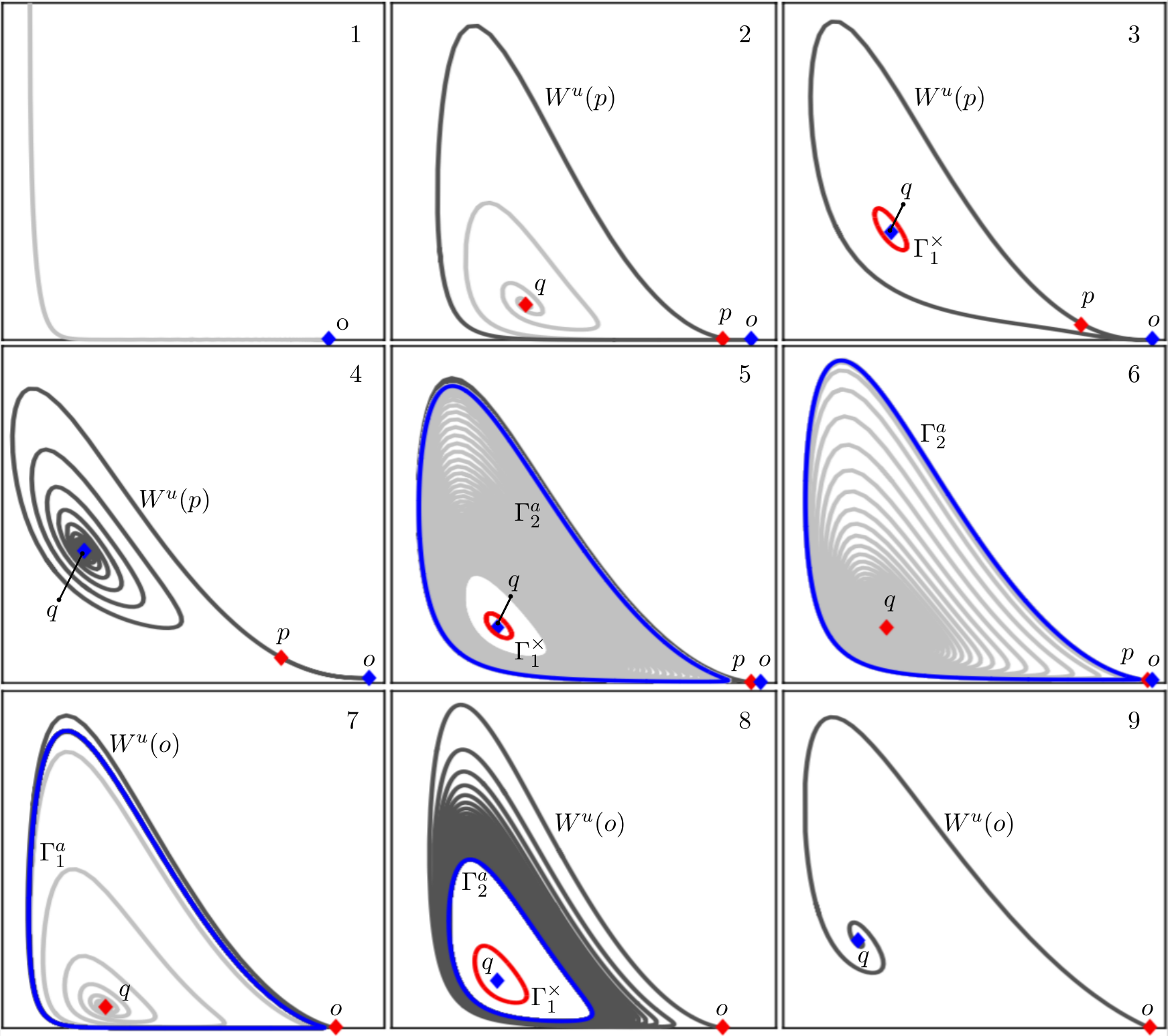}
\caption{The nine phase portraits projected onto the $(G, I)$-plane showing equilibria (blue and red diamonds) and periodic orbits (blue and red curves), one-dimensional manifolds (dark grey curves) and selected trajectories (light grey curves); attracting objects are blue and saddle objects are red. The chosen parameters are listed in Table~\ref{regions1to11}.}  
\label{Phase_portraits}
\end{figure}

\begin{figure}[t!]
\centering
\includegraphics[width=1 \columnwidth]{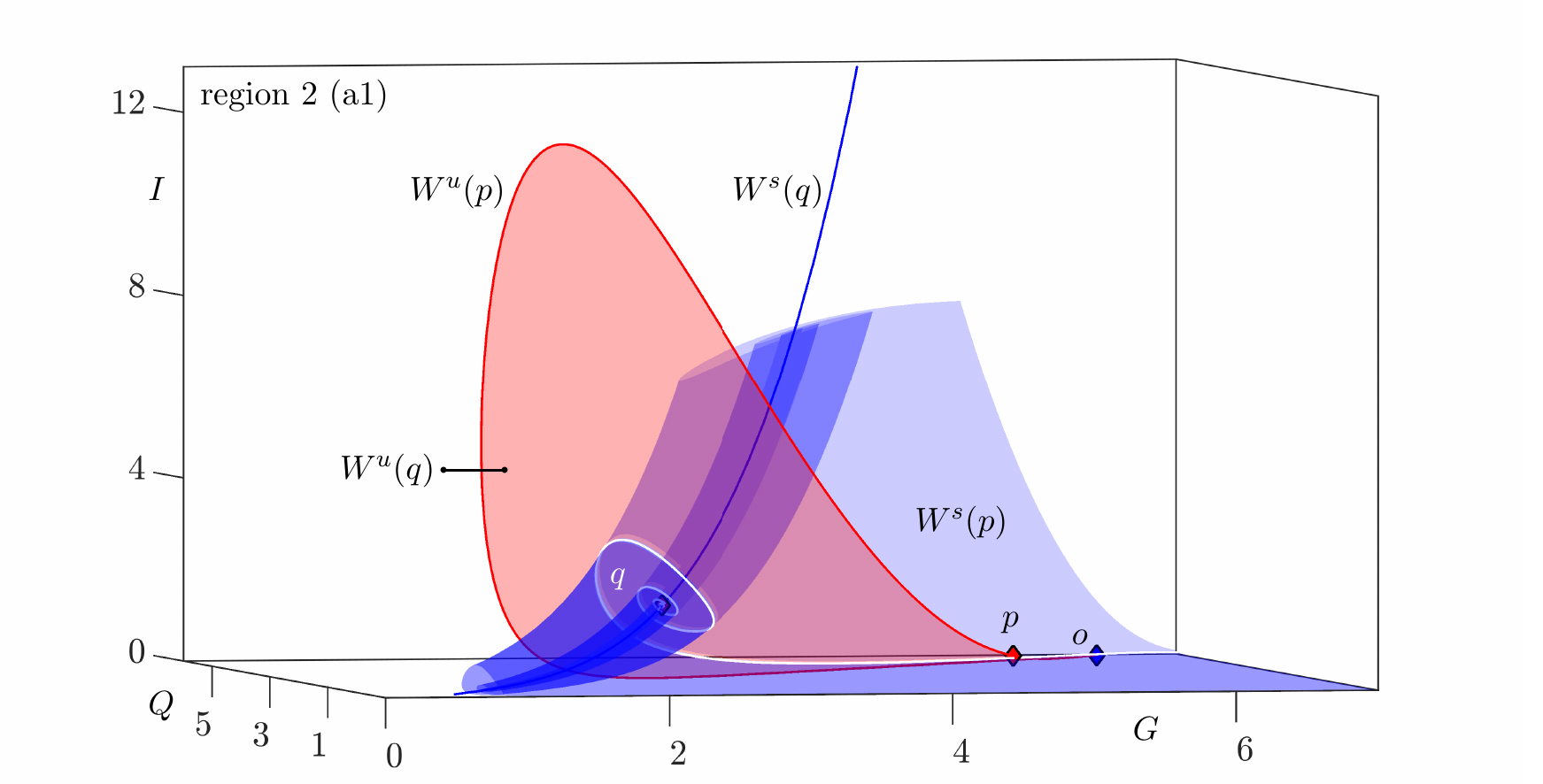}
\caption{Phase portrait in region 2 in the $(G,Q,I)$-space showing the attracting equilibrium $o$ (blue diamond), saddle equilibria $p$ and $q$ (red diamonds), manifolds $W^{u}(p)$ (red curve), $W^{s}(q)$ (blue curve), $W^{s}(p)$ (blue surface) and $W^{u}(q)$ (red surface); see Table~\ref{regions1to11} for parameter values.}
\label{3D_region2}
\end{figure}

\begin{table}[h]
\caption{Chosen representative parameter values for the phase portraits in the different regions of Fig.~\ref{Phase_portraits} and Fig.~\ref{stand_alone}, with an indication of invariant objects, namely equilibria $p, q, o$, and periodic orbits $\Gamma _{}^{a}$ (attractor) and $\Gamma _{}^{\times}$ (saddle); throughout, $B$=5.8 and $a$=1.8.}
\centering
{\begin{tabular}{cc|c|c|c|c|c|cl}\\[-2pt]
\toprule
\multicolumn{1}{c}{\multirow{2}[0]{*}{REGION }} & \multicolumn{1}{c}{\multirow{2}[0]{*}{($A$, $\gamma_{G}$, $\sigma$)}} 	 & \multicolumn{7}{c}{\kern-2em INVARIANT OBJECTS} 	\\[1pt]
								&       & $p$     & $q$     & $o$     & $\Gamma _{2}^{a}$    & $\Gamma _{1}^{\times}$     & $\Gamma _{3}^{\times}$   &   \\ 	\\[-7pt]
\hline																								\\

    1     & (5.500, 0.04000, 1.000)  &      		 &      		& \checkmark     &      		&      		&  \\
    2     & (6.540, 0.04000, 1.000)  & \checkmark    & \checkmark    & \checkmark     &      		&      		&  \\
    3     & (6.351, 0.12000, 1.000)  & \checkmark    & \checkmark    & \checkmark     &      		& \checkmark     &  \\
    4     & (6.720, 0.15000, 1.000)  & \checkmark    & \checkmark    & \checkmark     &      		&      		&  \\
    5     & (6.300, 0.07450, 1.000)  & \checkmark    & \checkmark    & \checkmark     & \checkmark    &      		& \checkmark  \\
    6     & (6.766, 0.06741, 1.000)  & \checkmark    & \checkmark    & \checkmark     & \checkmark    &      		&  \\
    7     & (6.800, 0.03000, 1.000)  &    			& \checkmark    & \checkmark     & \checkmark    &      		&  \\
    8     & (6.885, 0.06300, 1.000)  &      		& \checkmark    & \checkmark     & \checkmark    &      		& \checkmark \\
    9     & (7.500, 0.15000, 1.000)  &      		& \checkmark    & \checkmark     &      		&      		&  \\
    10    & (6.537, 0.04304, 1.115) & \checkmark    & \checkmark    & \checkmark     & \checkmark    & \checkmark    & \checkmark \\
    11   & (6.542, 0.04255, 1.115)  & \checkmark    & \checkmark    & \checkmark     & \checkmark    &      		& \checkmark \\
\end{tabular}}
  \label{regions1to11}%
\end{table}

\subsection{Phase portraits 1 to 9}
\label{sec:phase1to9}

The two-dimensional sketches of Fig.~\ref{Bifurcation_Sketches_original_study} complete the topological information needed for the bifurcation diagram. They have been confirmed in~\cite{Dubbeldam_Krauskopf_Self_pulsations_lasers_saturable_absorber} by computed versions in projection onto the $(G,I)$-plane consisting of equilibria, periodic orbits, one-dimensional invariant manifolds and selected trajectories. Figure~\ref{Phase_portraits} shows this type of planar representation, which suffices to characterize the topological aspects of the dynamics of the system. On the other hand, it leaves out considerable subtleties in the three-dimensional phase space, especially regarding how two-dimensional manifolds organize the phase space geometrically. 

\begin{figure}[t!]
\centering
\includegraphics[width=1 \columnwidth]{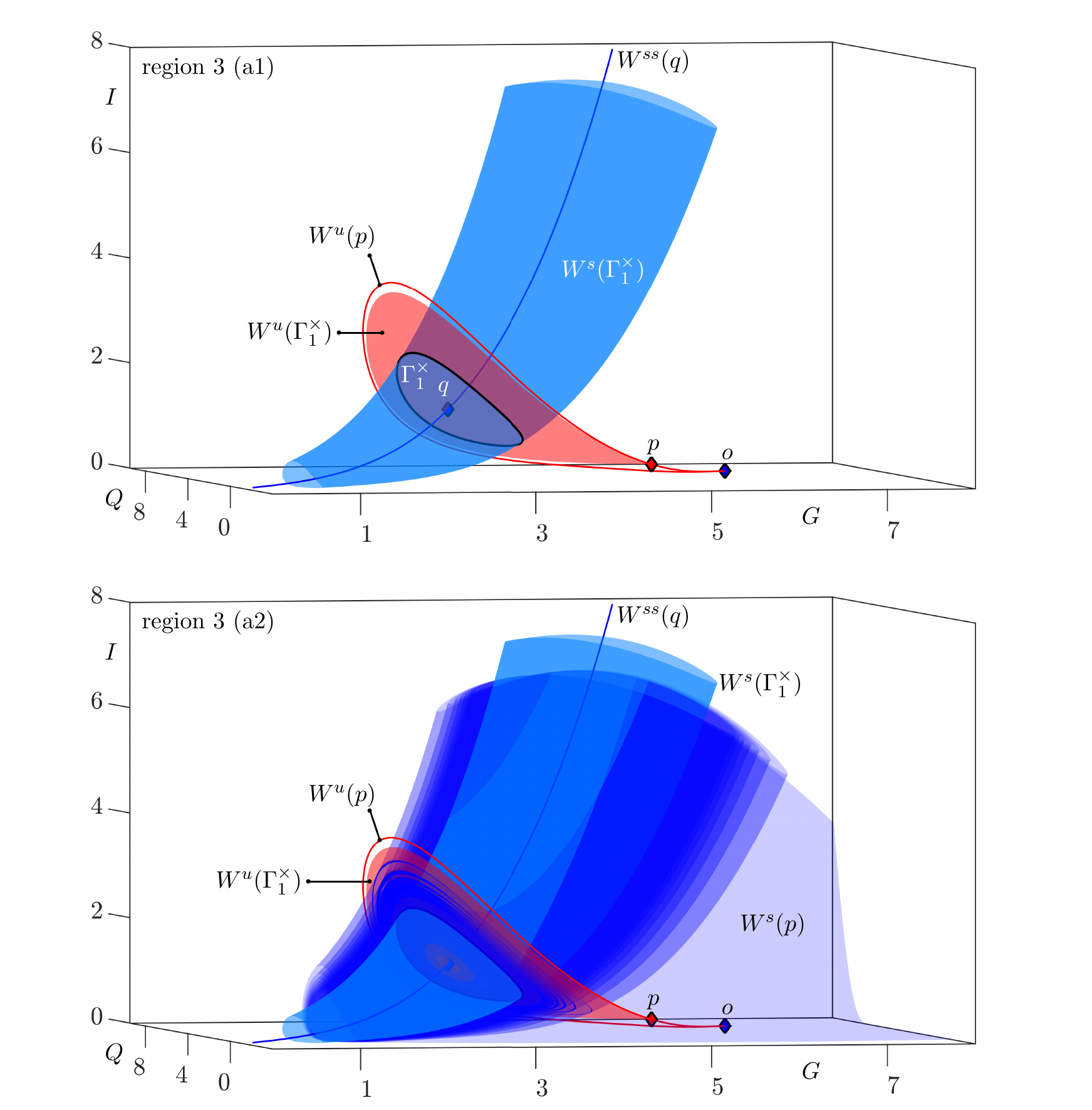}
\caption{Phase portrait in region 3 in the $(G,Q,I)$-space  showing the attracting equilibria $o$ and $q$ (blue diamond), saddle equilibrium $p$ (red diamond), saddle periodic orbit $\Gamma _{1}^{\times}$ (black curve) and manifolds $W^{ss}(q)$ (blue curve) and $W^{u}(p)$ (red curve). Panel (a1) shows the surfaces $W^{s}(\Gamma _{1}^{\times})$ (light blue surface) and $W^{u}(\Gamma _{1}^{\times})$ (red surface) and panel (a2) also shows $W^{s}(p)$ (dark blue surface); see Table~\ref{regions1to11} for parameter values.}  
\label{3D_region3}
\end{figure}

\begin{figure}[t!]
\centering
\includegraphics[width=1 \columnwidth]{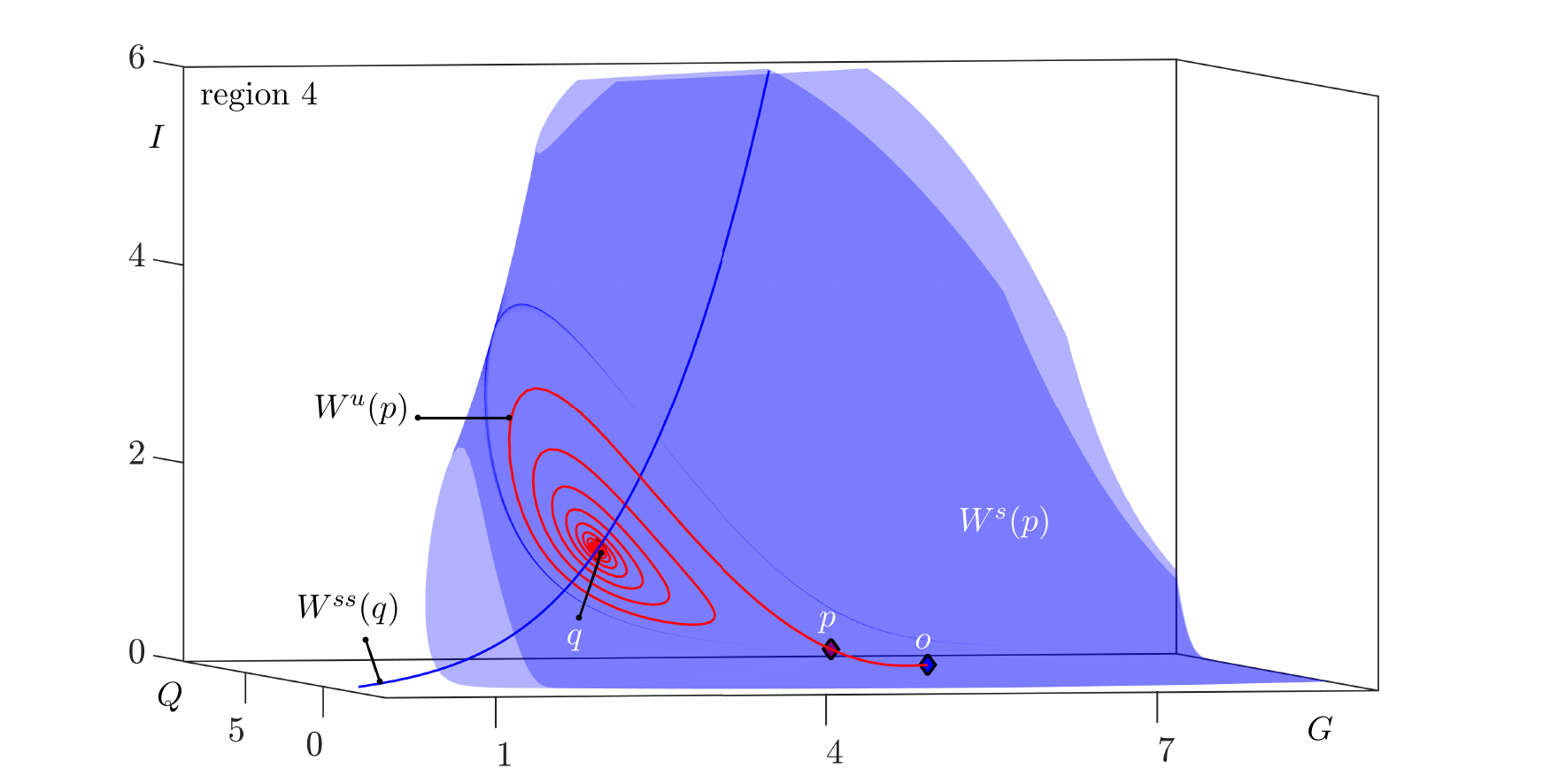}
\caption{Phase portraits in region 4 in the $(G,Q,I)$-space  showing the attracting equilibria $o$,$q$ (blue diamond), saddle equilibrium $p$ (red diamond) and manifold $W^{ss}(q)$ (blue curve). Also shown are the surfaces $W^{s}(p)$ (blue surface); see Table~\ref{regions1to11} for parameter values.}  
\label{3D_region4}
\end{figure}

\begin{figure}[t!]
\centering
\includegraphics[width=1 \columnwidth]{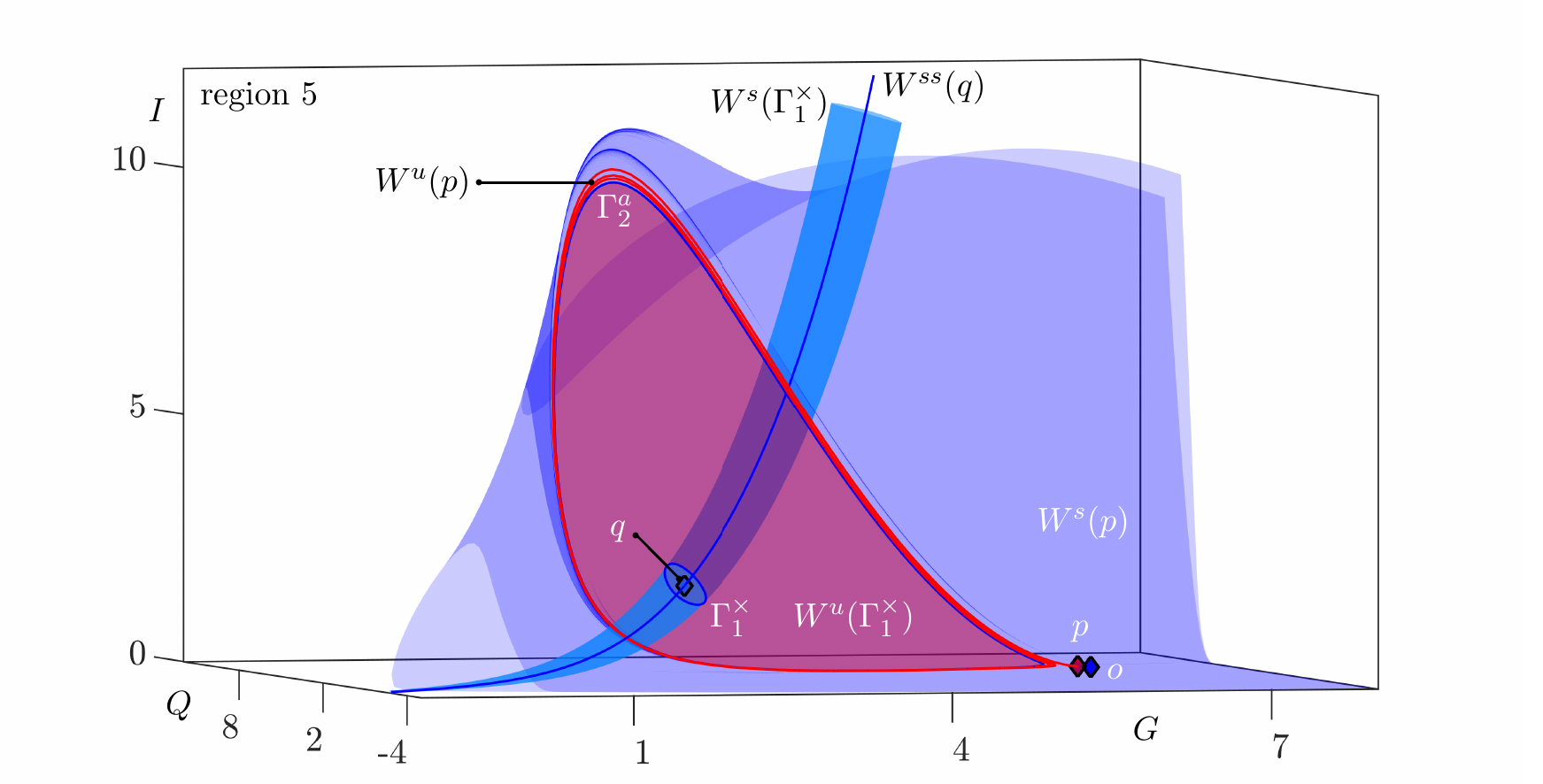}
\caption{Phase portrait in region 5 in the $(G,Q,I)$-space  showing the attracting equilibria $o$ and $q$ (blue diamond), saddle equilibrium $p$ (red diamond), attracting periodic orbit $\Gamma _{2}^{a}$ (blue curve), saddle periodic orbit $\Gamma _{1}^{\times}$ (black curve) and manifolds $W^{u}(p)$ (red curve) and $W^{ss}(q)$ (blue curve). Also shown are surfaces $W^{s}(p)$ (dark blue surface), $W^{s}(\Gamma _{1}^{\times})$ (light blue surface) and $W^{u}(\Gamma _{1}^{\times})$ (red surface); see Table~\ref{regions1to11} for parameter values.}  
\label{3D_region5}
\end{figure}

\begin{figure}[t!]
\centering
\includegraphics[width=1 \columnwidth]{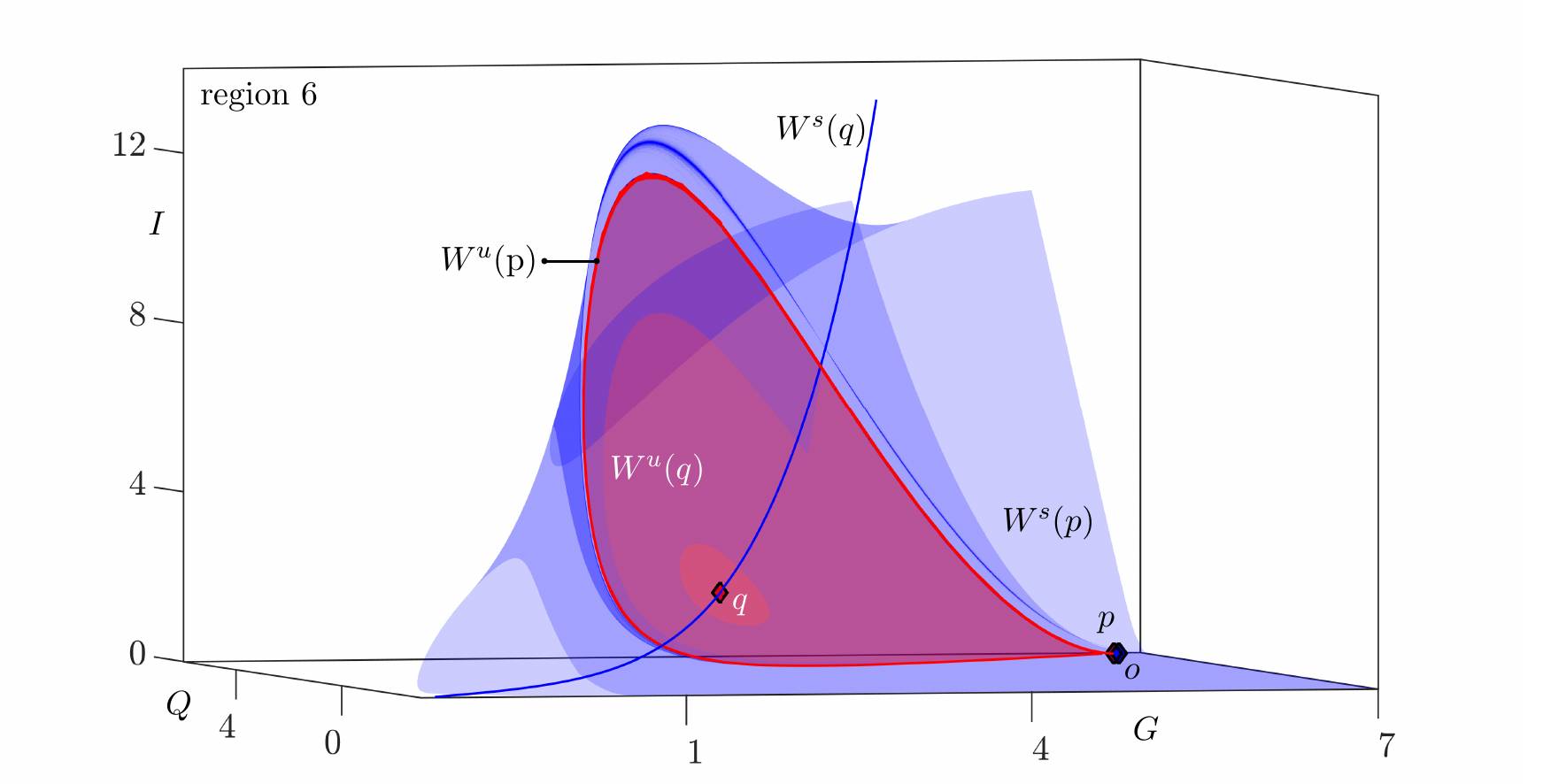}
\caption{Phase portrait in region 6 in the $(G,Q,I)$-space  showing the attracting equilibrium $o$ (blue diamond), saddle equilibria $q$ and $p$ (red diamond), attracting periodic orbit $\Gamma _{2}^{a}$ (blue curve) and manifolds $W^{u}(p)$ (red curve) and $W^{s}(q)$ (blue curve). Also shown are surfaces $W^{s}(p)$ (blue surface) and $W^{u}(q)$ (red surface); see Table~\ref{regions1to11} for parameter values.}  
\label{3D_region6}
\end{figure}

\begin{figure}[t!]
\centering
\includegraphics[width=1 \columnwidth]{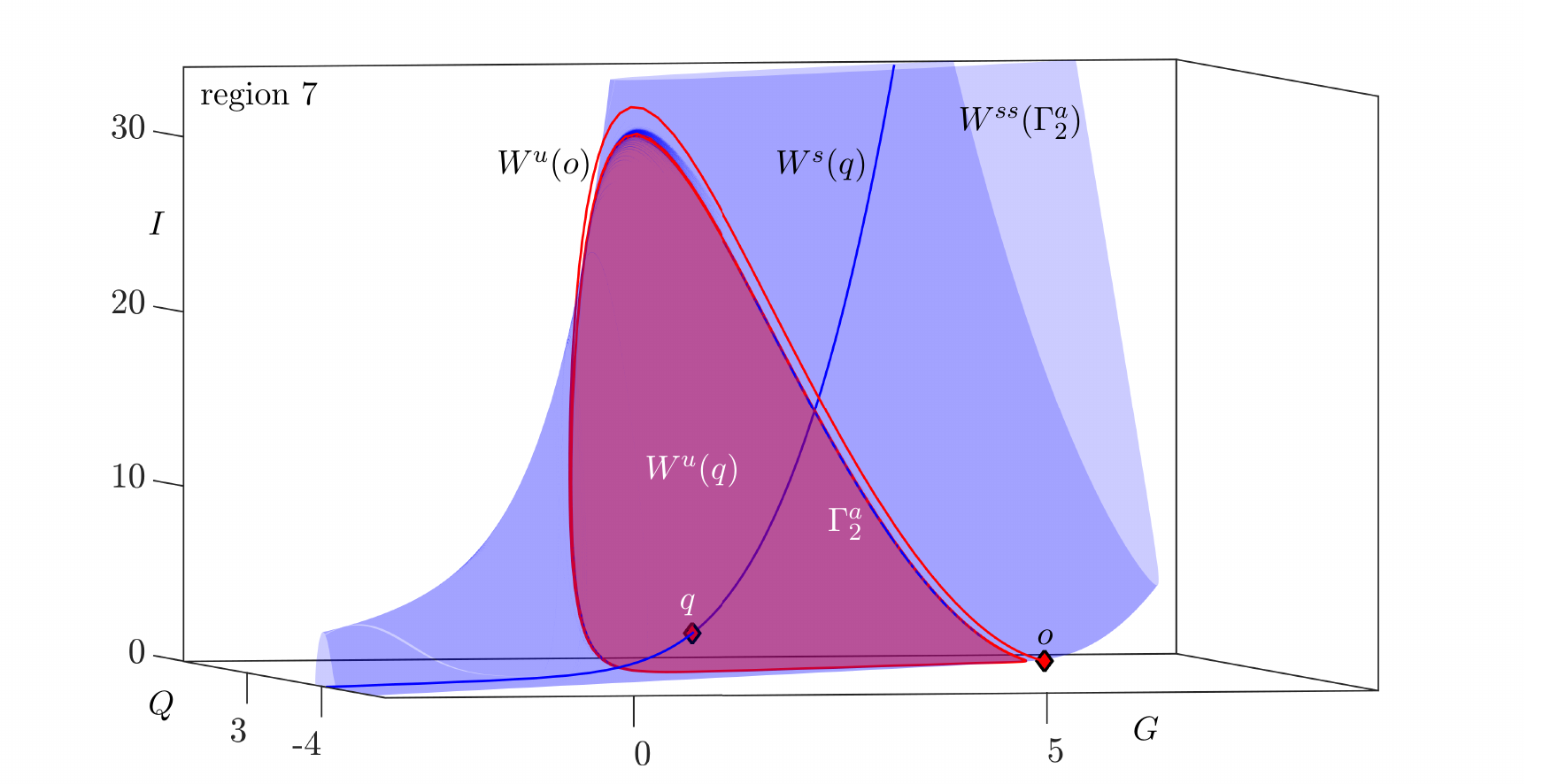}
\caption{Phase portrait in region 7 in the $(G,Q,I)$-space  showing the saddle equilibria $o$ and $p$ (red diamond), attracting periodic orbit $\Gamma _{2}^{a}$ (blue curve), saddle periodic orbit $\Gamma _{1}^{\times}$ (black curve) and manifolds $W^{u}(o)$ (red curve) and $W^{s}(q)$ (blue curve). Also shown are surfaces $W^{u}(q)$ (red surface) and $W^{ss}(\Gamma _{2}^{a})$ (blue surface); see Table~\ref{regions1to11} for parameter values.}  
\label{3D_region7}
\end{figure}

\begin{figure}[h!]
\centering
\includegraphics[width=1 \columnwidth]{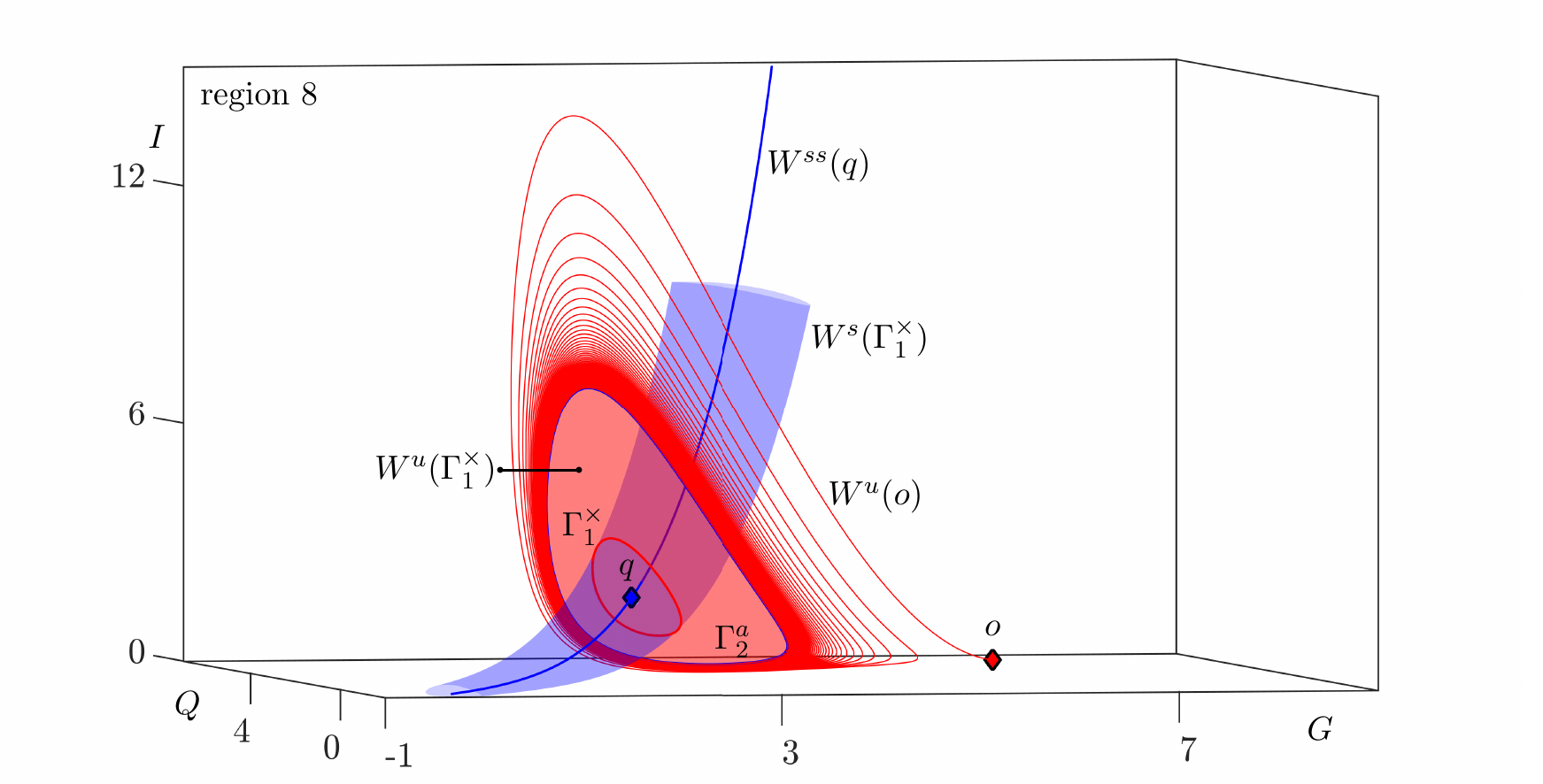}
\caption{Phase portrait in region 8 in the $(G,Q,I)$-space  showing the attracting equilibrium $q$ (blue diamond), saddle equilibrium $o$ (red diamond), attracting periodic orbit $\Gamma _{2}^{a}$ (blue curve), saddle periodic orbit $\Gamma _{1}^{\times}$ (red curve) and manifolds $W^{u}(o)$ (red curve), $W^{ss}(q)$ (blue curve). Also shown are surfaces $W^{s}(\Gamma _{1}^{\times})$ (blue surface) and $W^{u}(\Gamma _{1}^{\times})$ (red surface); see Table~\ref{regions1to11} for parameter values.}  
\label{3D_region8}
\end{figure}

\begin{figure}[t!]
\centering
\includegraphics[width=1 \columnwidth]{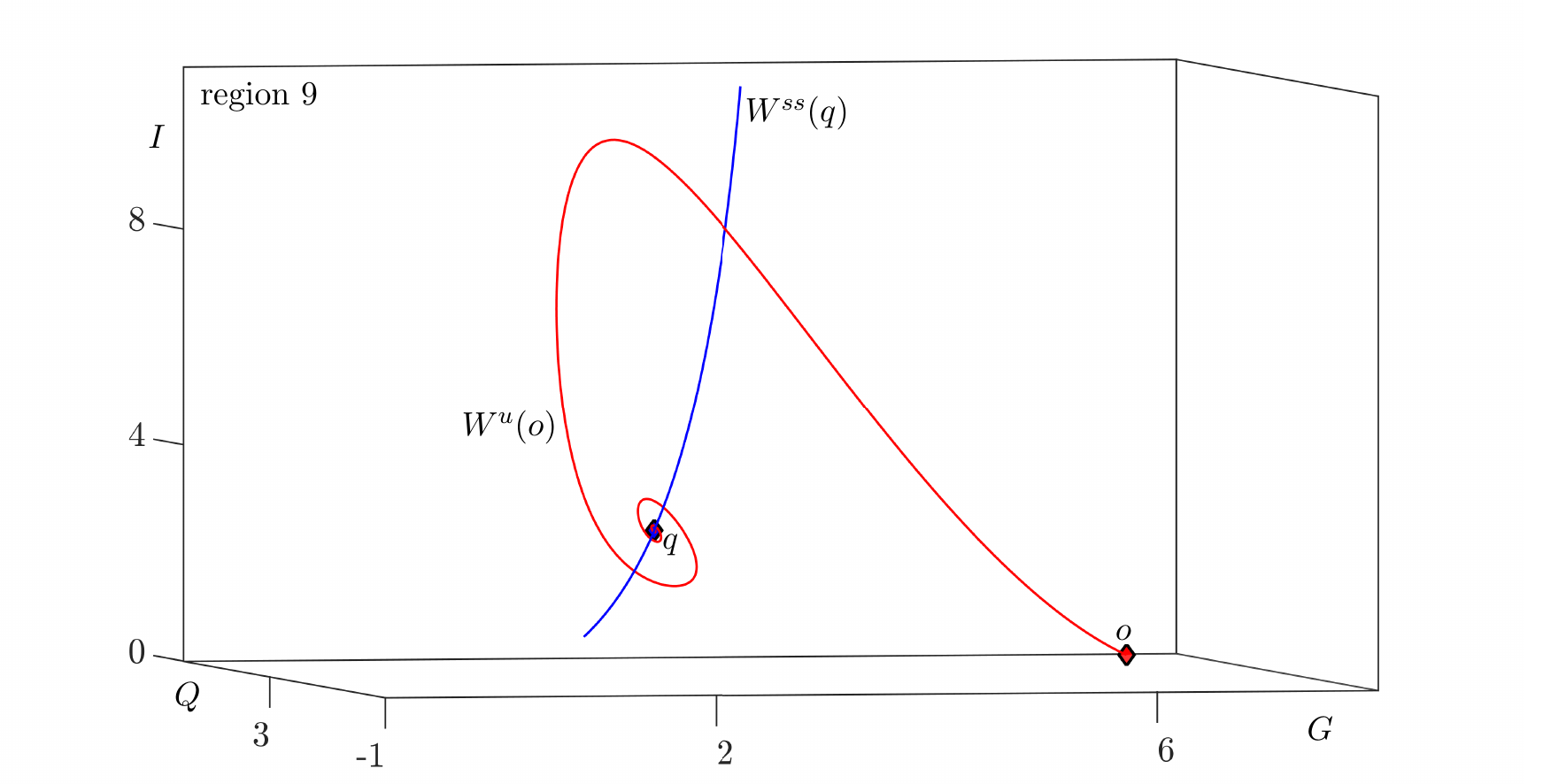}
\caption{Phase portraits in region 9 in the $(G,Q,I)$-space  showing the attracting equilibrium $q$ (blue diamond), saddle equilibrium $o$ (red diamond) and manifolds $W^{u}(o)$ (red curve), $W^{ss}(q)$ (blue curve); see Table~\ref{regions1to11} for parameter values.}  
\label{3D_region9}
\end{figure}

Figures~\ref{3D_region2} to \ref{3D_region9} show phase portraits 2 to 9 in the full $(G, Q, I)$-space with all relevant one-dimensional and two-dimensional invariant manifolds; all phase portraits shown have been computed for parameter values detailed in Table~\ref{regions1to11}.
We proceed by discussing the nine generic phase portraits 1 to 9 in detail, with reference to the sketches in Fig.~\ref{Bifurcation_Sketches_original_study}, the two-dimensional computed phase portraits in Fig.~\ref{Phase_portraits} and their three-dimensional computed equivalents in $(G, Q, I)$-space in figures~\ref{3D_region2} to \ref{3D_region9}. 

To the left of the curve S is region 1, where the stable equilibrium $o$ on the invariant plane $I=0$ is the only equilibrium; compare with panel 1 of Fig.~\ref{Phase_portraits} (because of its simplicity we are not presenting a three-dimensional phase portrait for region 1). The point $o$ represents the off-state of the laser. In region 1, this equilibrium is stable and attracts all initial conditions.

Crossing the curve S below the point BT means moving into region 2. There are now two additional equilibria $p$ and $q$, where $p$ is a saddle with two stable and one unstable eigenvalues and $q$ is a spiral source. The off-state $o$ is still the only attractor and attracts all initial conditions except for $q$, $p$ and its stable manifold $W^{s}(p)$. The sketch in Fig.~\ref{Bifurcation_Sketches_original_study} as well as panel 2 of Fig.~\ref{Phase_portraits} show the equilibria and the one-dimensional unstable manifold $W^{u}(p)$ of $p$, each side of which ends at the attractor $o$. Also shown in panel 2 of  Fig.~\ref{Phase_portraits} is a trajectory starting near $q$; it spirals out and also ends up at $o$, which is indeed the only attractor. Since point $p$ is of saddle type it must have a two-dimensional stable manifold $W^{s}(p)$. It is indicated in the sketch of Fig.~\ref{Phase_portraits} as two trajectories ending up at $p$. As Fig.~\ref{3D_region2} shows, $W^{s}(p)$ is a surface that rolls up in a carpet-like fashion around the one-dimensional stable manifold $W^{s}(q)$ of $q$. Again, $W^{s}(q)$ is not shown in panel 2 of Fig.~\ref{Phase_portraits}; it can be thought of here and in the sketch of Fig.~\ref{Bifurcation_Sketches_original_study} as the third, stable direction that that is not shown in these two-dimensional representations. Importantly, the surface $W^{s}(p)$ forms a separatrix that defines the excitability threshold of the system in region 2. Any initial condition near point $p$ and below $W^{s}(p)$ goes back to the off state $o$ immediately, whereas any initial condition above $W^{s}(p)$ makes a large excursion along $W^{u}(p)$  before ending up at $o$. This all-or-nothing response is known as excitability and the laser is excitable in region 2. Note further that $W^{s}(p)$ intersects the red two-dimensional unstable manifold $W^{u}(q)$ of $q$ to form a structurally stable heteroclinic connection, which is shown as a white curve in Fig.~\ref{3D_region2}.

The transition from region 2 to 3 is via the curve H. The points $o$ and $p$ along with its manifolds $W^{u}(p)$ and $W^{s}(p)$  remain largely unchanged. However, the point $q$ changes from saddle type to an attractor with the appearance of a periodic orbit, $\Gamma^{\times}_{1}$ of saddle type around $q$; see panel 3 of Fig.~\ref{Bifurcation_Sketches_original_study} and~\ref{Phase_portraits}. The equilibrium $q$ possesses a one-dimensional strong stable manifold $W^{ss}(q)$ which again corresponds to the third stable direction that has been left out in this projection. There are now two attractors in region two. Since the periodic orbit $\Gamma^{\times}_{1}$ is of saddle type it has two-dimensional stable and unstable manifolds $W^{s}(\Gamma _{1}^{\times})$ and $W^{u}(\Gamma _{1}^{\times})$, respectively. They are shown in Fig.~\ref{3D_region3}(a1) as blue and red surfaces. In fact, $W^{s}(\Gamma _{1}^{\times})$ forms the basin boundary separating the two attractors $q$ and $p$. Panel (a2) of Fig.~\ref{3D_region3} shows also the stable manifold $W^{s}(p)$ which, similarly rolls up into a carpet-like fashion, but accumulates now on the stable manifold $W^{s}(\Gamma _{1}^{\times})$ of the periodic orbit $\Gamma^{\times}_{1}$.

Crossing the homoclinic bifurcation L from region 3 into region 4 results in the disappearance of the periodic orbit $\Gamma^{\times}_{1}$. The sketches of Fig.~\ref{Bifurcation_Sketches_original_study} show that in region 4, the equilibria $o$, $p$ and $q$ keep their previous stability properties but the structure of the one-dimensional manifold $W^{u}(p)$ has changed; see also panel 4 of Fig.~\ref{Phase_portraits}, where the equilibria $o$ and $q$ remain the only attractors. One branch of $W^{u}(p)$ goes to the attractor $o$ while the other branch of $W^{u}(p)$ spirals towards the attractor $q$. The two-dimensional manifold $W^{s}(p)$ of the saddle point $p$ has been sketched as a curve in Fig.~\ref{Bifurcation_Sketches_original_study} due to the chosen projection; this curve corresponds to the basin boundary between the two attractors $o$ and $q$ and the laser is bistable in region 4. The actual nature of the basin $W^{s}(p)$ boundary is illustrated in Fig.~\ref{3D_region4}. The stable manifold $W^{s}(p)$ folds over to envelope the point $q$ and the spiralling part of $W^{u}(p)$ that ends up at $q$, while the shorter part of $W^{u}(p)$ that goes to the off-state $o$ remains outside and below $W^{s}(p)$. As such, any initial condition on the q-side of $W^{s}(p)$ goes to the point $q$ and the laser produces constant intensity output, while any initial condition on the other side of $W^{s}(p)$ ends up at $o$, that is, at the off-state with zero intensity.

Moving into region 5 requires crossing the curve D. The difference now is the appearance of two periodic orbits: an attracting periodic orbit $\Gamma^{a}_{2}$ and the periodic orbit $\Gamma^{\times}_{1}$ of saddle type; see panels 5 of Fig.~\ref{Bifurcation_Sketches_original_study} and~\ref{Phase_portraits}. This brings to three the total number of attractors and the unstable manifold $W^{u}(p)$ now accumulates onto the outside of the periodic orbit $\Gamma^{a}_{2}$ from the outside; see also the trajectory in panel 5 of Fig~\ref{Phase_portraits} in region 5 starting near $\Gamma^{\times}_{1}$ and spiralling out to end up on $\Gamma^{a}_{2}$. Notice that the two periodic orbits are far apart and $\Gamma^{a}_{2}$ is close to homoclinic. The presence of these three attractors requires basin boundaries that are not comprehensively shown in figures~\ref{Bifurcation_Sketches_original_study} and~\ref{Phase_portraits} but rather in the three dimensional view of Fig.~\ref{3D_region5}. Here we see that the surfaces $W^{s}(p)$ and $W^{s}(\Gamma _{1}^{\times})$  form two different basin boundaries that delineate three basins and, as a consequence, the laser is tri-stable. The manifold $W^{s}(\Gamma _{1}^{\times})$ forms a cylinder and any initial condition inside it goes to q, while any initial condition on the outside of $W^{s}(\Gamma _{1}^{\times})$ but on the inner side of $W^{s}(p)$ goes to $\Gamma^{a}_{2}$ which corresponds to a periodic, pulsing laser intensity. Additionally, all initial conditions on the other-side of $W^{s}(p)$ go to $o$.

Region 6 is then entered via the curve H. The difference here is the disappearance of the periodic orbit $\Gamma^{\times}_{1}$  at the curve H; the equilibrium $q$ is now of saddle type as sketch 6 of Fig~\ref{Bifurcation_Sketches_original_study} illustrates. The two-dimensional manifold $W^{s}(p)$, which is sketched as a curve due to the chosen projection, defines a separatrix between the two attractors, $\Gamma^{a}_{2}$ and $o$. Hence the laser is bistable in region 6 and, thus, capable of pulsations or being off, depending on the initial condition. Any initial condition close to the equilibrium $q$ spirals out towards the attractor $\Gamma^{a}_{2}$ as represented by the light grey curve in panel 6 of Fig~\ref{Phase_portraits}. Figure~\ref{3D_region6} shows the geometry of the one-dimensional stable manifold $W^{s}(q)$ and the separatrix $W^{s}(p)$ in the full phase space. The manifold $W^{s}(p)$ forms a loop around, $\Gamma^{a}_{2}$ and, consequently, any initial condition on the inside of $W^{s}(p)$ goes to $\Gamma^{a}_{2}$, while any initial condition on the other-side of $W^{s}(p)$ goes to the off-state $o$.

Crossing the curve T from region 6 leads into region 7. At T the point $p$ disappears and the attractor $o$ now becomes a saddle in region 7. Hence, $o$ now has a two-dimensional stable manifold $W_{s}(o)$ which actually is the plane defined by $I=0$. Its one-dimensional unstable manifold $W^{u}(p)$ accumulates onto the periodic orbit $\Gamma^{a}_{2}$ as is sketched in Fig~\ref{Bifurcation_Sketches_original_study}; see also panel 7 of Fig.~\ref{Phase_portraits}. Similarly, any initial condition close to the equilibrium $q$ spirals out quickly towards the periodic orbit $\Gamma^{a}_{2}$, which is the only attractor in region 7. Figure~\ref{3D_region7} shows geometry in $\mathbb{R}^3$.  Also shown is the two dimensional strong stable manifold $W^{ss}(\Gamma _{2}^{a})$ around which  $W^{u}(o)$ approaches $\Gamma^{a}_{2}$.

Crossing the curve H from region 7 leads to region 8. The saddle point $q$ is now an attractor and the periodic orbit $\Gamma^{\times}_{1}$ of saddle type now sits within $\Gamma^{a}_{2}$ and around $q$. In the sketch of region 8 in Fig.~\ref{Bifurcation_Sketches_original_study} and in panel 8 of Fig.~\ref{Phase_portraits}, the periodic orbit  $\Gamma^{\times}_{1}$ separates the two only attractors in region 8 where we again find bistability. The two dimensional manifold $\Gamma^{\times}_{1}$ is not shown here but it delineates the two basins of attraction $q$ and $\Gamma^{a}_{2}$. Figure~\ref{3D_region8} illustrates $W^{s}(\Gamma _{1}^{\times})$ in $\mathbb{R}^3$  as a cylindrical separatrix such that any initial condition inside $W^{s}(\Gamma _{1}^{\times})$ goes to $q$, while any initial condition on the outside of $W^{s}(\Gamma _{1}^{\times})$ goes to $\Gamma^{a}_{2}$; the latter includes $W^{u}(o)$.

Finally, moving into region 9 from region 7 is also via the curve H.  The sketch of region 9 in Fig.~\ref{Bifurcation_Sketches_original_study} and panel 9 of Fig.~\ref{Phase_portraits} show that the periodic orbit $\Gamma^{a}_{2}$ has disappeared and $q$ is now stable and the only attractor. Hence all initial conditions end up at $q$ and the laser produces constant lasing output. Figure~\ref{3D_region9} shows   this relatively simple phase portrait in the three-dimensional phase space.

The bifurcation diagram in Fig.~\ref{Bifurcation_Sketches_original_study} is structurally stable, so changing $\sigma$ only a little bit causes no change in its topological structure. Larger changes beyond a certain $\sigma$, on the other hand, may results in topological changes of the bifurcation diagram. For decreasing $\sigma$ from $\sigma = 1$, which in physical terms means that the gain decays slower than the absorption, we found in our investigation that the bifurcation diagram does not change qualitatively. In other words, for $\sigma \in (0, 1]$ the relative positions of the bifurcation curves remain the same, the different regions change only in size and none of them disappears. This means that bifurcation diagram case BI as shown in Fig.~\ref{Bifurcation_Sketches_original_study} is relevant whenever $\gamma_{G} \leq \gamma_{Q}$.

\section{Bifurcation diagrams BII to BX for $1 < \sigma$ }

We now increase $\sigma$ from 1 so that the gain section always decays faster than the absorber section, that is, $\gamma_{Q} < \gamma_{G}$. There are now qualitative changes to the relative arrangement of bifurcation curves and new dynamics is created.  Every such change of the bifurcation diagram corresponds to an event of codimension three. These come in two flavours: codimension-two-plus-one or codimension-three. The first case usually involves a codimension-two point of the bifurcation diagram moving across a codimension-one curve; however, there is no codimension-three bifurcation in the sense that the two bifurcations do not involve the same object in phase space. A codimension-three bifurcation in the classical sense of bifurcation theory, on the other hand, concerns a single object in phase space and typically generates additional bifurcation curves. We will see several examples, including the codimension-three event where the BT point occurs at $I=0$, defining the invariant plane. In either case, a codimension-three event creates changes that occur locally in some smaller area of the bifurcation diagram; in other words, outside this neighbourhood the bifurcation diagram remains unaltered qualitatively.

\begin{table}[!htb]
\caption{$\sigma$-values for the shown bifurcation cases BI to BX, with locations of codimension-two points in Fig.~\ref{Bif_diagram_sigma}; throughout, $B$=5.8 and $a$=1.8.}
\centering
{ \begin{tabular}{c c|c|c|c|c}\\[-2pt]
\toprule
	CASE 	&$\sigma$	&BT ($A$, $\gamma_G$,) 			&C ($A$, $\gamma_G$,) 	 		&NH ($A$, $\gamma_G$,) 				&GH ($A$, $\gamma_G$,)		\\[1pt]  
\hline																																			\\[-2pt]
       BI     		& 1.0000		& (6.1300, 0.2100)				& n/a 	            		& (6.6080, 0.0880) 							& (7.2561, 0.0460)			\\[2pt]
       BII     		& 1.0670  		& (6.2118, 0.1366)				& (6.5869, 0.0648) 	            		& (6.6079, 0.0621) 							& (6.8149, 0.0462)			\\[2pt]
       BIII     	& 1.0697  		& (6.2176, 0.1334) 				& (6.5843, 0.0638) 	            		& (6.6078, 0.0608)  							& (6.7821, 0.0473) 			\\[2pt]
       BIV     	& 1.0700  		& (6.2183, 0.1331) 				& (6.5840, 0.0636) 	           		& (6.6078, 0.0606)  							& (6.7821, 0.0471) 			\\[2pt]
       BV     		& 1.1150  		& (6.3119, 0.0788) 				& (6.5307, 0.0437) 	& \multicolumn{1}{p{7.0em}|}{(6.6242, 0.0322) (6.7924, 0.0124)} 	&  (6.6001, 0.0375) 			\\[2pt]
       BV--VI     	& 1.121807  		& (6.3255, 0.0703) 				& (6.5210, 0.0400) 	            		& (6.6543, 0.0230) 							& (6.5753, 0.0351)  			\\[2pt]
       BVI     	& 1.1490  		& (6.3784, 0.0351) 				& (6.4791, 0.0225) 	            		& n/a 							& (6.4902, 0.0216) 			\\[2pt]	
       BVI--VII     & 1.149328  		& (6.3789, 0.0348) 				& (6.4786, 0.0223) 	           		& n/a 							& (6.4893, 0.0213) 			\\[2pt]	
       BVII     	& 1.1534  		& (6.3802, 0.0339) 				& (6.4768, 0.0219) 	          		& n/a 							& (6.4875, 0.0209) 			\\[2pt]	
       BVIII     	& 1.1550  		& (6.3896, 0.0273) 				& (6.4760, 0.0180) 	         			& n/a 							& (6.4754, 0.0175) 			\\[2pt]
       BVIII--BIX & 1.1754978	& (6.4274, 3.561$\times10^{-7}$)	& (6.4274, 3.561$\times10^{-7}$) 	& n/a 							& (6.4274, 3.561$\times10^{-7}$)\\[2pt]
       BIX	     	& 1.1755  		& n/a				& n/a 	         			& n/a 							& n/a 			\\[2pt]	
       BIX--BX    	& 1.5353044  		& n/a				& n/a 	         			& n/a 							& n/a 			\\[2pt]
       BX 	   	& 1.550  		& n/a				& n/a 	         			& n/a 							& n/a 			\\[2pt]
\end{tabular}}
 \label{bif_cases}%
\end{table}

\begin{figure}[t!]
\centering
\includegraphics[width=1\columnwidth]{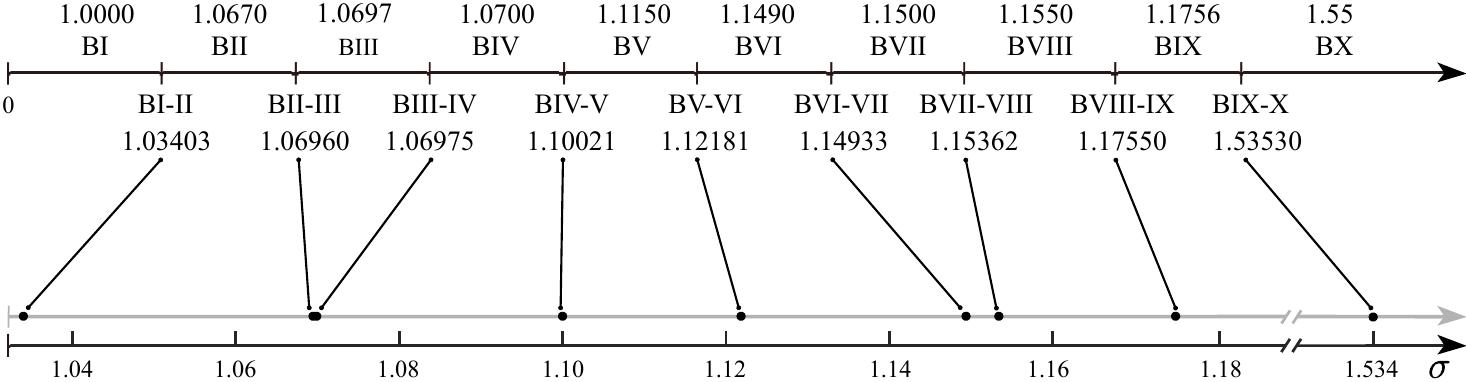}
\caption{A sketch of the $\sigma$-line (top) indicating the relative positions of all bifurcation cases BI--BX with values chosen for the computed images in Fig.~\ref{Bif_diagram_sigma}; identified values of the associated transitions are listed and indicated on the actual $\sigma$-line (bottom).}
\label{sigma_values_bif_types}
\end{figure}

Apart from bifurcation diagram BI, there are nine other bifurcation diagrams, which we denote BII up to BX. Our results were derived with topological arguments from bifurcation theory in combination with a careful numerical investigation. Each transition between successive bifurcation diagrams is a codimension-three event that we will identify and present as part of our discussion. In subsequent sections we present topological sketches of cases BII to BX, as well as computed bifurcation diagrams in the $(A,\gamma_G)$-plane for values of the parameter $\sigma$ listed in Table~\ref{bif_cases}. As an overview and outlook at what follows, Fig.~\ref{sigma_values_bif_types} shows the $\sigma$-line with all cases BI--BX and the associated transitions between them, both as a sketch at the top and with actual positions of numerically computed transitions at the bottom. 

\begin{figure}[t!]
\centering
\includegraphics[width=0.9\columnwidth]{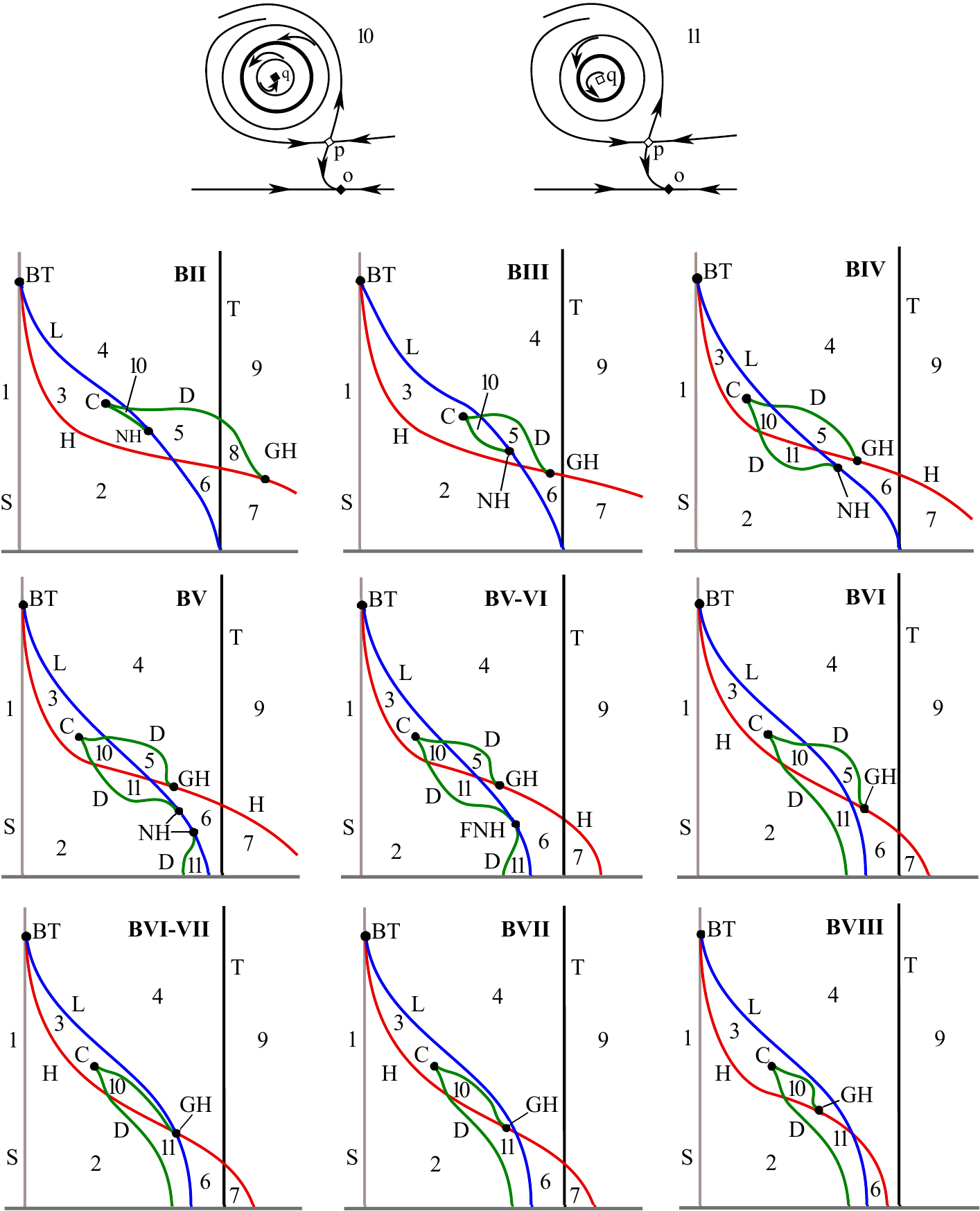}
\caption{Sketches of bifurcation diagrams in the $(A,\gamma_{G})$-plane of generic cases \textbf{BII} to \textbf{BVIII}, and transitions \textbf{BV--VI} and \textbf{BVI--VII}, with phase portraits for the two new regions 10 and 11; compare with Fig.~\ref{Bifurcation_Sketches_original_study}.}
\label{Bifurcation_Sketches}
\end{figure}

\begin{figure}[t!]
\centering
\includegraphics[width=1 \columnwidth]{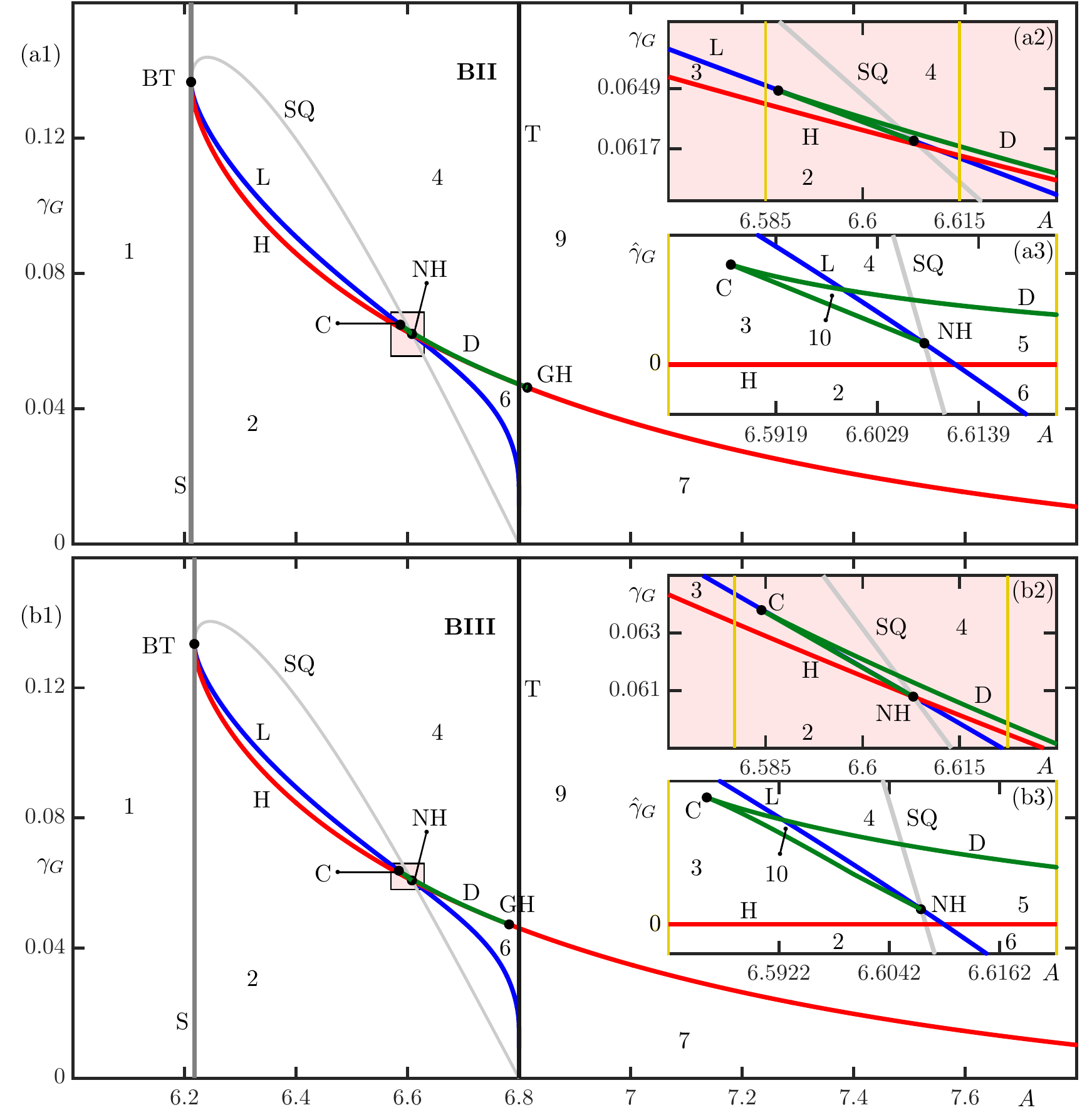}
\caption{Computed bifurcation diagram in the $(A,\gamma_{G})$-plane of case BII at $\sigma$=1.0670 in panel (a) and for case BIII at $\sigma=1.0697$ in panel (b). Shown are curves S of saddle-nodes of equilibria (dark grey), T of transcritical (black), D of saddle-node of periodic orbits (green), H of Hopf bifurcations (red) and L of homoclinic bifurcations (blue). Along the curve SQ (light grey) the saddle quantity is zero. Also shown are points of codimension-two bifurcations, namely points of generalised Hopf (GH), neutral-saddle homoclinic (NH), Cusp (C) and Bogdanov-Takens (BT) bifurcations. Inserts (a2) and (b2) show enlargements of the respective shaded regions, and inserts (a3) and (b3) show the bifurcation curves plotted relative to curve H in the region bounded by the vertical yellow lines in (a2) and (b2), respectively.}  
\label{Bif_diagram_sigma}
\end{figure}

Cases BII to BVIII, which all feature a Bogdanaov-Takens point and associated bifurcation curves, are discussed in detail in Sec.~\ref{sec:ItoVIII}, where we present topological sketches in Fig.~\ref{Bifurcation_Sketches} and associated numerical evidence in Fig.~\ref{Bif_diagram_sigma}; also shown here are two-dimensional and three-dimensional images of the new phase portraits 10 and 11 in Fig.~\ref{stand_alone} and Figs.~\ref{3D_region10} and ~\ref{3D_region11}, repectively. The final two cases BIX and BX, which do no longer feature a Bogdanaov-Takens point, are discussed in Sec.~\ref{sec:IXtoX}, where we explain them with a projection of the relevant bifurcation loci onto the $(A,\sigma)$-plane.

\renewcommand{\thefigure}{\arabic{figure} (continued)}
\addtocounter{figure}{-1}
\begin{figure}[t!]
\centering
\includegraphics[width=1 \columnwidth]{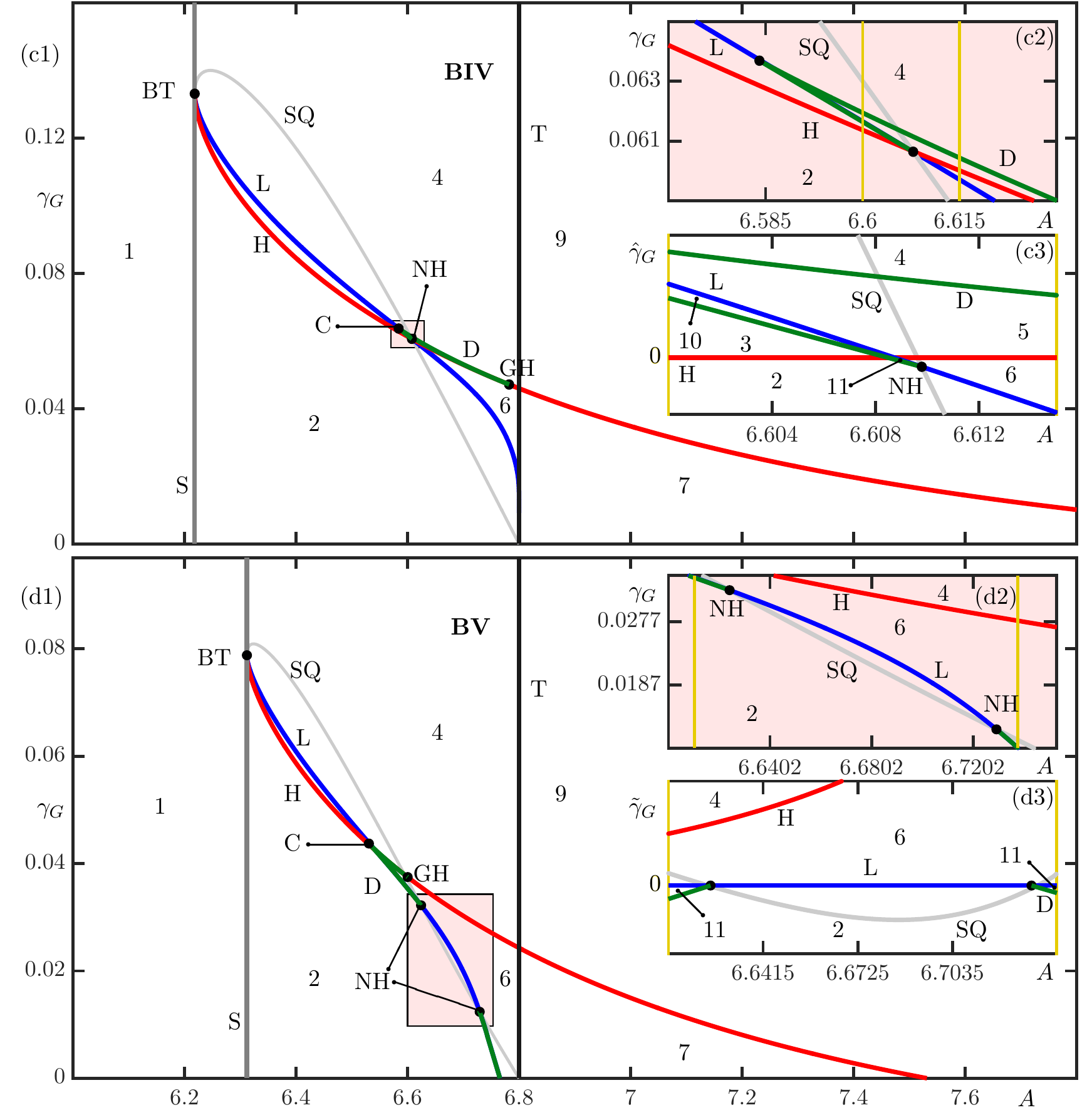}
\caption{Computed bifurcation diagram in the $(A,\gamma_{G})$-plane for case BIV at $\sigma=1.0700$ in panels (c) and for case BV at $\sigma=1.1150$ in panels (d). Here insert (c3) shows curves relative to H, while insert (d3) shows curves relative to L.}  
\end{figure}
\renewcommand{\thefigure}{\arabic{figure}}

\renewcommand{\thefigure}{\arabic{figure} (continued)}
\addtocounter{figure}{-1}
\begin{figure}[t!]
\centering
\includegraphics[width=1 \columnwidth]{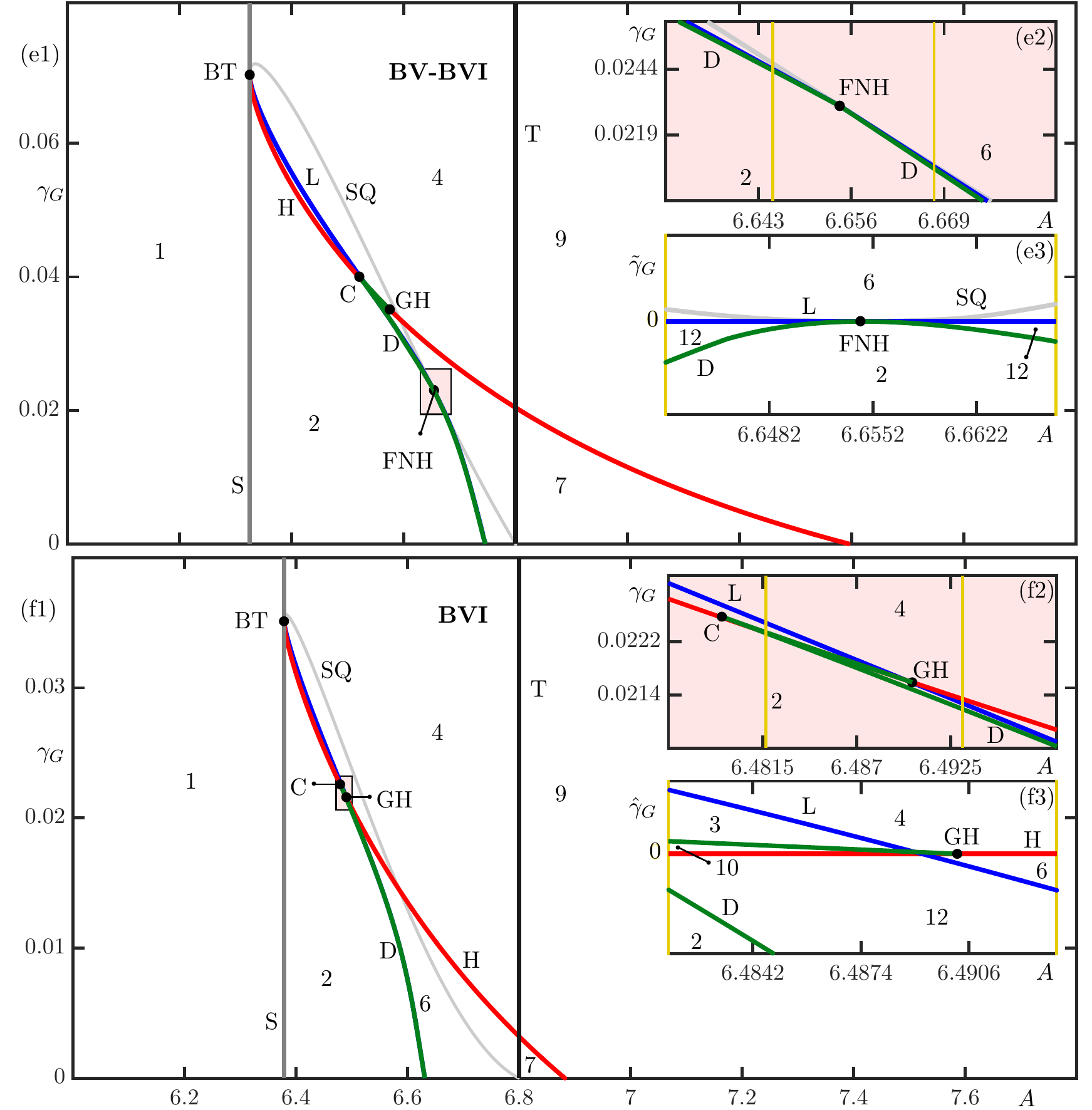}
\caption{Computed bifurcation diagram in the $(A,\gamma_{G})$-plane case BV--VI at $\sigma=1.12181$ in panels (e) and for case BVI at $\sigma=1.1490$ in panels (f). Here insert (e3) shows curves relative to L, while insert (f3) shows curves relative to H.}  
\end{figure}
\renewcommand{\thefigure}{\arabic{figure}}

\renewcommand{\thefigure}{\arabic{figure} (continued)}
\addtocounter{figure}{-1}
\begin{figure}[t!]
\centering
\includegraphics[width=1 \columnwidth]{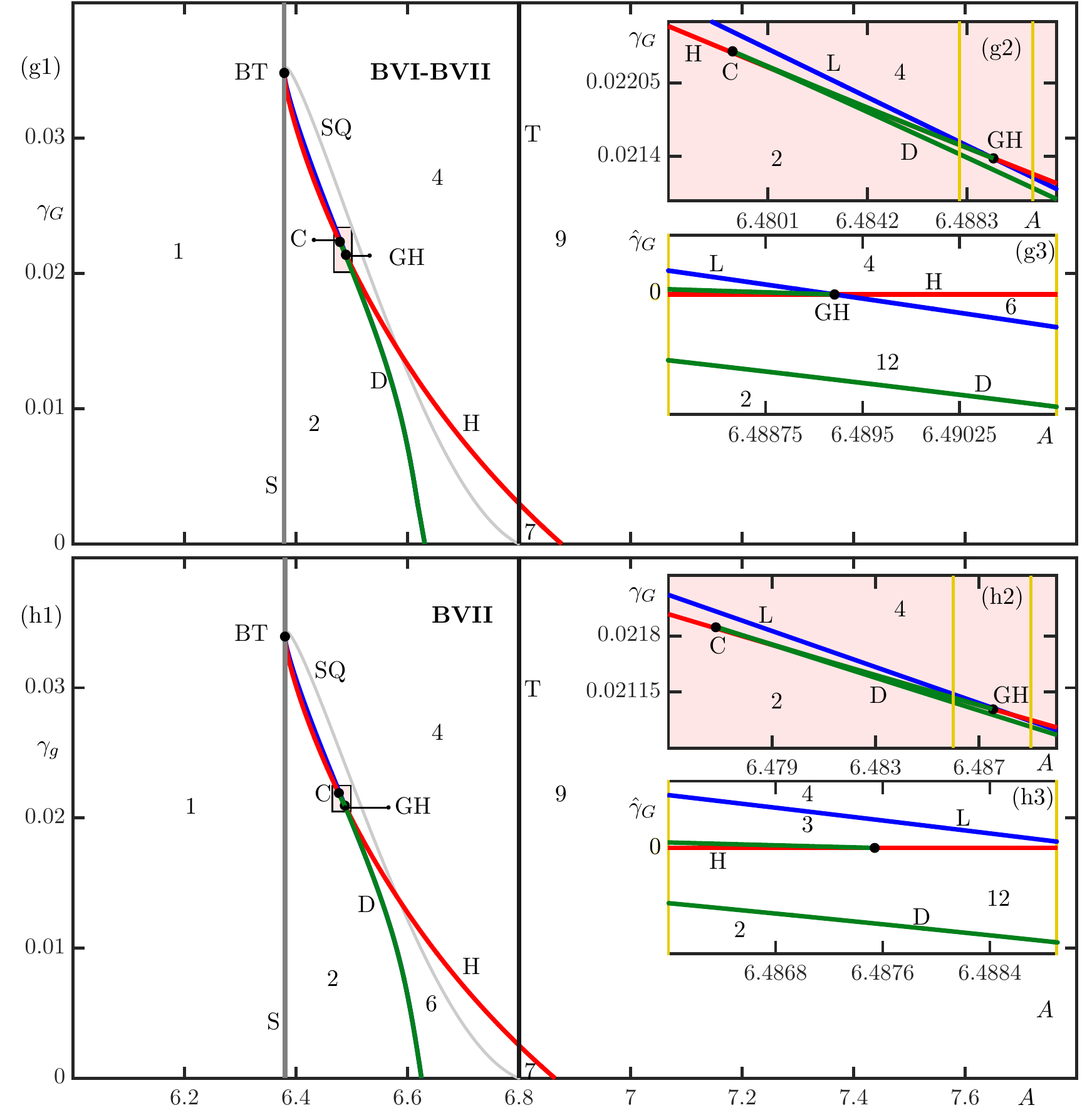}
\caption{Computed bifurcation diagram in the $(A,\gamma_{G})$-plane of transition BVI--VII at $\sigma=1.14933$ in panels (g) and for case BVII at $\sigma=1.1500$ in panels (h).}  
\end{figure}
\renewcommand{\thefigure}{\arabic{figure}}

\renewcommand{\thefigure}{\arabic{figure} (continued)}
\addtocounter{figure}{-1}
\begin{figure}[t!]
\centering
\includegraphics[width=0.95 \columnwidth]{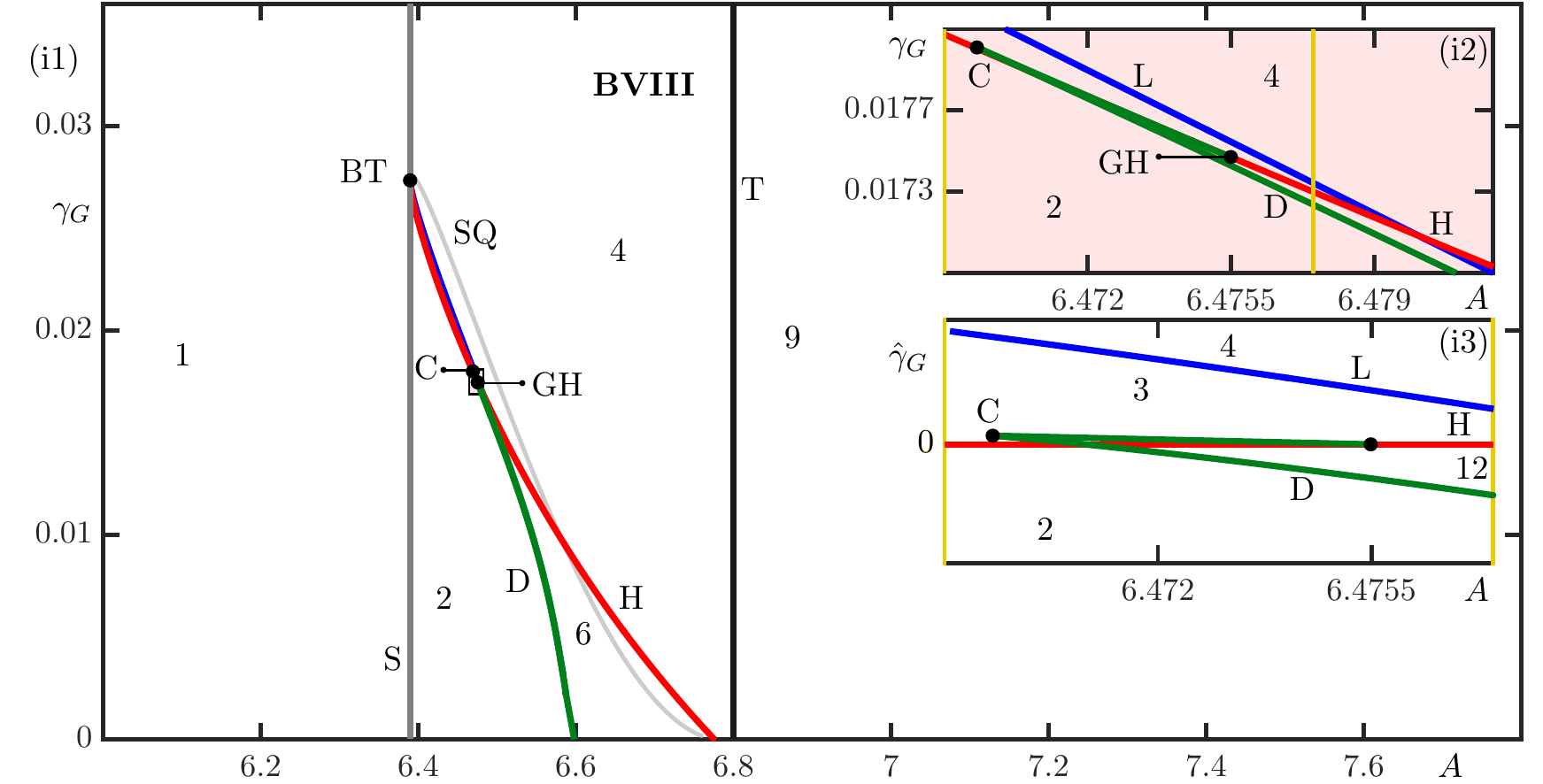}
\caption{Computed bifurcation diagram in the $(A,\gamma_{G})$-plane of case BVIII at $\sigma=1.1550$. }  
\end{figure}
\renewcommand{\thefigure}{\arabic{figure}}

\subsection{Cases BII to BVIII}
\label{sec:ItoVIII}

As $\sigma$ is increased from $\sigma = 1$, the first codimension-three event we encounter is the change in type of the point NH at a neutral-saddle homoclinic bifurcation (where the saddle quantity is zero). This occurs at $\sigma = 1.03403$ when the part of the curve D tangent to curve L switches from the top-right-hand side of L to the bottom-left-hand side of L. This forces the appearance of a point C of cusp bifurcation on the curve D of saddle-node bifurcations of periodic orbits, so that D can remain connected to the rest of the bifurcation diagram above it. This transition gives the new case BII shown in Fig~\ref{Bifurcation_Sketches}. The bifurcation diagram is unchanged except for the curve D now crossing  the curve L, forming a cusp and connecting to the point NH from the other-side. As a result, the new region 10 is created. To show that this is really the case, Fig.~\ref{Bif_diagram_sigma}(a) at $\sigma=1.0670$ presents the numerical evidence for the bifurcation diagram of case BII. Panel (a1) provides a view of the entire bifurcation diagram. The different bifurcation curves are topologically as in the sketch in Fig.~\ref{Bifurcation_Sketches}, but lie very close to one another near the points C and NH. Also shown  is the curve SQ. It is not actually a bifurcation curve  (for that reason it is shown in light grey) and it does not delineate any regions; however where it crosses the curve L is exactly  where the point NH occurs from which the curve D emerges. Panel (a2) shows an enlargement near the points C and NH, enabling us see better the various regions. The curves H, L, and D are a bit clearer, but still close to one another. To clarify their positioning, panel (a3) shows the same curves plotted relative to the curve H. This is achieved by showing the signed distance $\hat{\gamma}_G$ in the $\gamma_G$-direction between the respective curves and the curve H. Hence, H is the $A$-axis in the $(\hat{\gamma}_G, A)$-plane of panel (a3), and regions 5, 6 and the newly formed region 10 are now visible very clearly. This shows that for bifurcation case BII it is possible to transition from region 3 to region 10 via the curve D. The sketch of the phase portrait in region 10 is shown at the top of Fig.~\ref{Bifurcation_Sketches}. Compared to region 3, two new periodic orbits exist in region 10: the attracting periodic orbit $\Gamma^{a}_{2}$ and the periodic orbit $\Gamma^{\times}_{3}$ of saddle type. This brings the total number of attractors to three. Phase portrait 10 will be discussed in more detail in Sec.~\ref{sec:phase10to11} with the help of computed two-dimensional and three-dimensional images.

As $\sigma$ is increased further, the next topological change is that the point GH moves from the right to the left of the curve T. This is a codimension-two-plus-one event at which region 8 of case BII disappears; however, no other change of the bifurcation diagram takes place at $\sigma=1.06960$. The bifurcation diagram of case BIII is like case BII with the only exception that the point GH is now completely to the left of T and region 8 has disappeared. The numerical evidence shown in Fig.~\ref{Bif_diagram_sigma}(b1) clearly illustrates that GH has moved from the right to the left of curve T. The enlargements (a2) and (b2) show that all the other curves are still positioned as for case BII. Panel (b3) shows the curves H, L and D plotted relative to the curve H in the $(A,\hat{\gamma}_G)$-plane.

The next transition as $\sigma$ is increased occurs when the point NH moves below the curve H. Again, it is a codimension-two-plus-one event when NH lies on the curve L and subsequently moves below H. This transition takes place at $\sigma=1.06975$. The bifurcation diagram of case BIV is shown in Fig~\ref{Bifurcation_Sketches} and Fig.~\ref{Bif_diagram_sigma}(c) presents the numerical evidence; its panel (c1) shows all bifurcation curves, which lie close to one another near the points C and NH. The enlargement of this region is shown in panel (c2): the curves H, L and D are still very close to one another but panel (c3) of the $(A,\hat{\gamma}_G)$-plane shows clearly that the point NH is now below the curve H. Notice closely regions 2, 5, 6, 10 and 11; the new region 11 is very small but the figure clearly shows that it indeed exists. Moving from region 10 to region 11 is achieved by crossing the curve H. The sketch of the phase portrait in region 11 at the top of  Fig~\ref{Bifurcation_Sketches} shows that the saddle periodic orbit $\Gamma^{\times}_{1}$ has disappeared and $q$ is now a saddle point. How this manifests itself in the three-dimensional phase space will be discussed in Sec.~\ref{sec:phase10to11}.

Increasing $\sigma$ further, we encounter another codimension-two-plus-one event when a second NH point emerges from $\gamma_{G}=0$. Consequently, a second region 11 is created as the NH point moves up along the curve L; see the sketch BV of Fig.~\ref{Bifurcation_Sketches}. Numerically, the transition from case BIV to case BV occurs at  $\sigma = 1.10021$, as determined by identifying the additional intersection point between L and SQ with $\gamma_G = 0$ as accurately as possible.
Figure~\ref{Bif_diagram_sigma}(d) presents the numerical evidence for the bifurcation diagram of case BV. Panel (d1) shows all bifurcation curves and, in particular, the two points NH. Indeed, the curve SQ intersects the curve L at these two points. Panel (d2) shows an enlargement near the two points NH, enabling us see better the regions delineated by the the curves H, L, and two curves D. Notice that the two curves D effectively lie on top of the curve L, even at this enlarged scale. This is why we plot in Fig.~\ref{Bif_diagram_sigma}(d3) the curves H, D and SQ relative to L by showing them in the $(A,\tilde{\gamma}_G)$-plane, where $\tilde{\gamma}_G$ is now the signed distance in the $\gamma_G$-direction between the respective curve and L (which is now the $A$-axis). Regions 2, 4, 6 and even the two regions 11 are now visible much more clearly.

The next transition occurs when these two neutral-saddle homoclinic points NH come together to form a tangency between the curves D and L. We denote this bifurcation a fold neutral-saddle homoclinic or FNH. The moment of this transition is shown in sketch BV--VI of Fig.~\ref{Bifurcation_Sketches}. We found numerically that this transition happens for $\sigma=1.12181$. The numerical evidence showing the transitional case BV--BVI is provided in Fig.~\ref{Bif_diagram_sigma}(e). The curves L, D and SQ are in very close proximity, with D and L practically indistinguishable in panel (e1); we remark that this remains to be so up to case BVIII. Even in the enlargement in panel (e2) near the point FNH shows the curves L, D and SQ in very close to one anothe, but the plot in Fig.~\ref{Bif_diagram_sigma}(e3) in the $(A,\tilde{\gamma}_G)$-plane illustrates clearly how all three curves are tangent at FNH. Therefore, for larger $\sigma$, case BVI of the bifurcation diagram emerges, which is sketched in Fig.~\ref{Bifurcation_Sketches} with the numerical evidence in  Fig~\ref{Bif_diagram_sigma}(f). The only topological difference from case BV is the disappearance of the two points NH and the merger of two regions 11 into a single region 11. Note that the point GH is now very close to the curve L, as shown in the $(A,\hat{\gamma}_G)$-plane of Fig.~\ref{Bif_diagram_sigma}(f3), resulting in a much smaller region 5. Finally, GH moves from the right to the left of L in a codimension-two-plus-one bifurcations shown as BVI--BVII in Fig.~\ref{Bifurcation_Sketches}. Numerically, we found that the point GH lies on the curve L for $\sigma=1.14933$. The numerical evidence for this transition is shown in Fig.~\ref{Bif_diagram_sigma}(g); see, in particular, the enlargements around the point GH in panels (g2) and (g3).

The resulting bifurcation diagram of case BVII for larger $\sigma$ is very similar to case BVI, with the only exception that the point GH is now to the left of the curve L. As a consequence, region 5 has disappeared; see also the numerical evidence in Fig.~\ref{Bif_diagram_sigma}(h). Panel (h2) shows that the point GH has completely moved from the right to the left of the curve L; see also panel (h3). The next transition we encounter, as $\sigma$ is increased further, occurs when the end point of the curve H on $\gamma_{G} = 0$ lies on the curve T. This is a codimension-two-plus-one event, which we refer to as HT$_0$. Past this special point, the end point of H is on the other side of the curve T, which results in the disappearance of region 7. We found numerically that HT$_0$ occurs at  $\sigma = 1.15362$. The bifurcation diagram of case BVIII is then as sketched in Fig.~\ref{Bifurcation_Sketches}, and the numerical evidence for it is presented in Fig.~\ref{Bif_diagram_sigma}(i). In panel (i1) the entire curve H is seen as lying to the left of T; the enlargements near the points GH and C in panels (i2) and (i3) show that there is no change otherwise.

\begin{figure}[t!]
\centering
\includegraphics[width=0.6 \columnwidth]{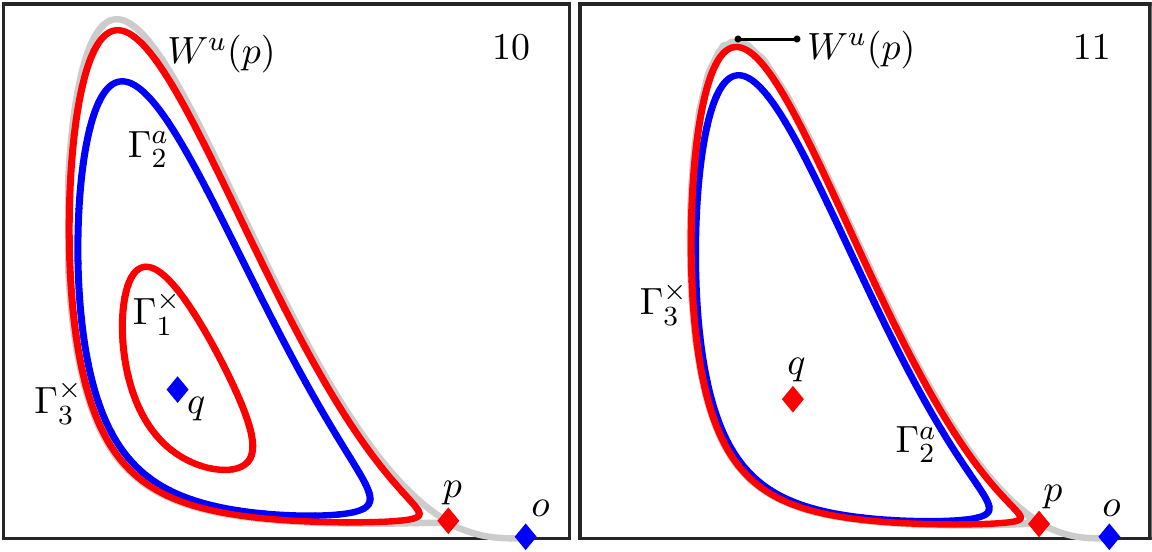}
\caption{The two phase portraits in regions 10 and 11 projected onto the $(G, I)$-plane, showing equilibria (blue and red diamonds) and periodic orbits (blue and red curves); attracting objects are blue and saddle objects are red.}
\label{stand_alone}
\end{figure}

\begin{figure}[t!]
\centering
\includegraphics[width=1 \columnwidth]{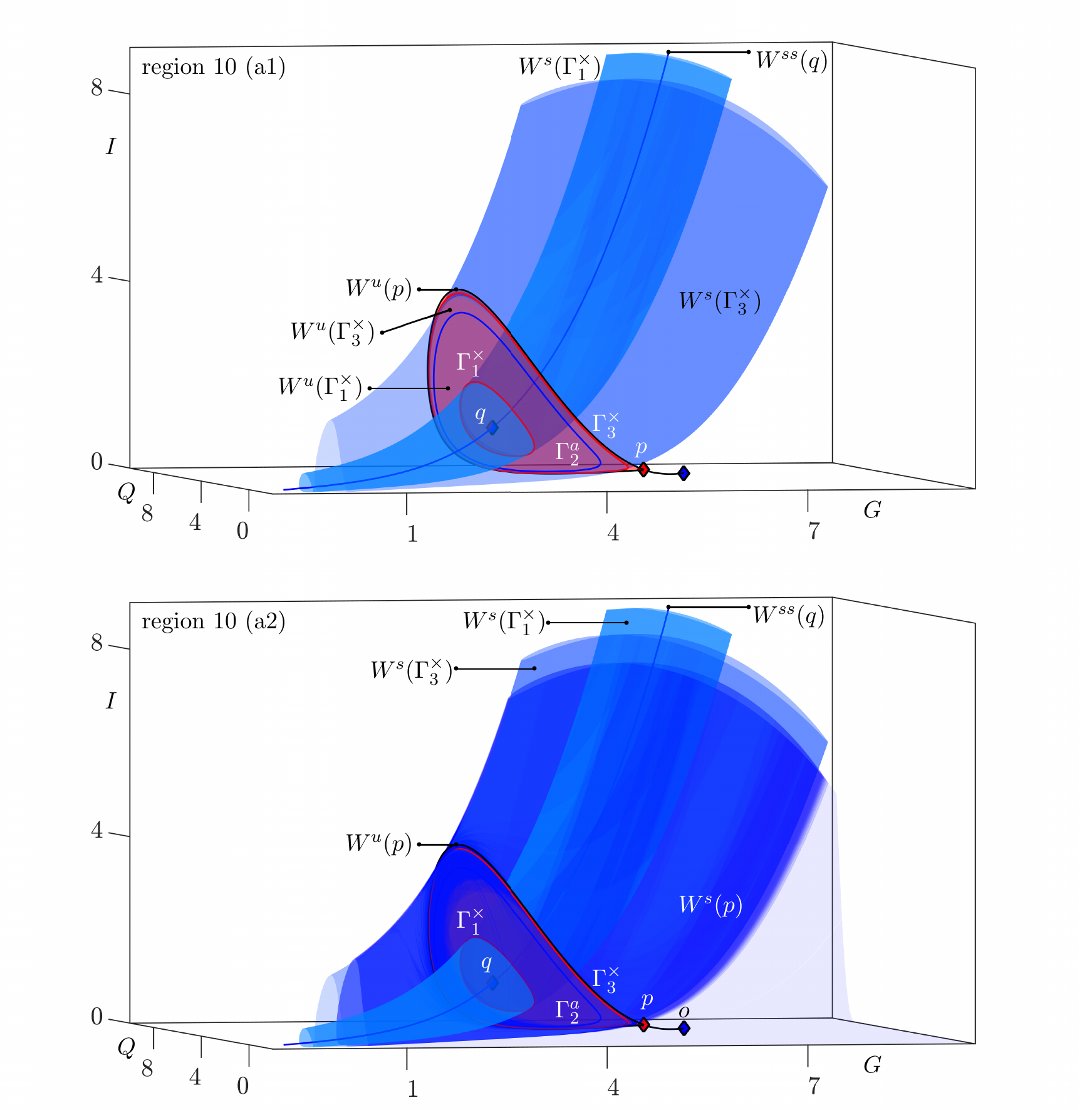}
\caption{Phase portrait in region 10 in the $(G,Q,I)$-space  showing the attracting equilibria $o$ and $q$ (blue diamond), saddle equilibrium $p$ (red diamond), attracting periodic orbit $\Gamma _{2}^{a}$ (blue curve), saddle periodic orbits $\Gamma _{1}^{\times}$ (red curve) and $\Gamma _{3}^{\times}$ (red curve) and manifolds $W^{ss}(q)$ (blue curve) and $W^{u}(o)$ (black curve). Panel (a1) shows the surfaces $W^{s}(\Gamma _{1}^{\times})$ (teal surface), $W^{s}(\Gamma _{3}^{\times})$ (light blue surface), $W^{u}(\Gamma _{1}^{\times})$ (red surface) and $W^{u}(\Gamma _{3}^{\times})$ (red surface), panel (a2) also shows $W^{s}(p)$ (dark blue surface); see Table~\ref{regions1to11} for parameter values.}
\label{3D_region10}
\end{figure}

\begin{figure}[t!]
\centering
\includegraphics[width=1 \columnwidth]{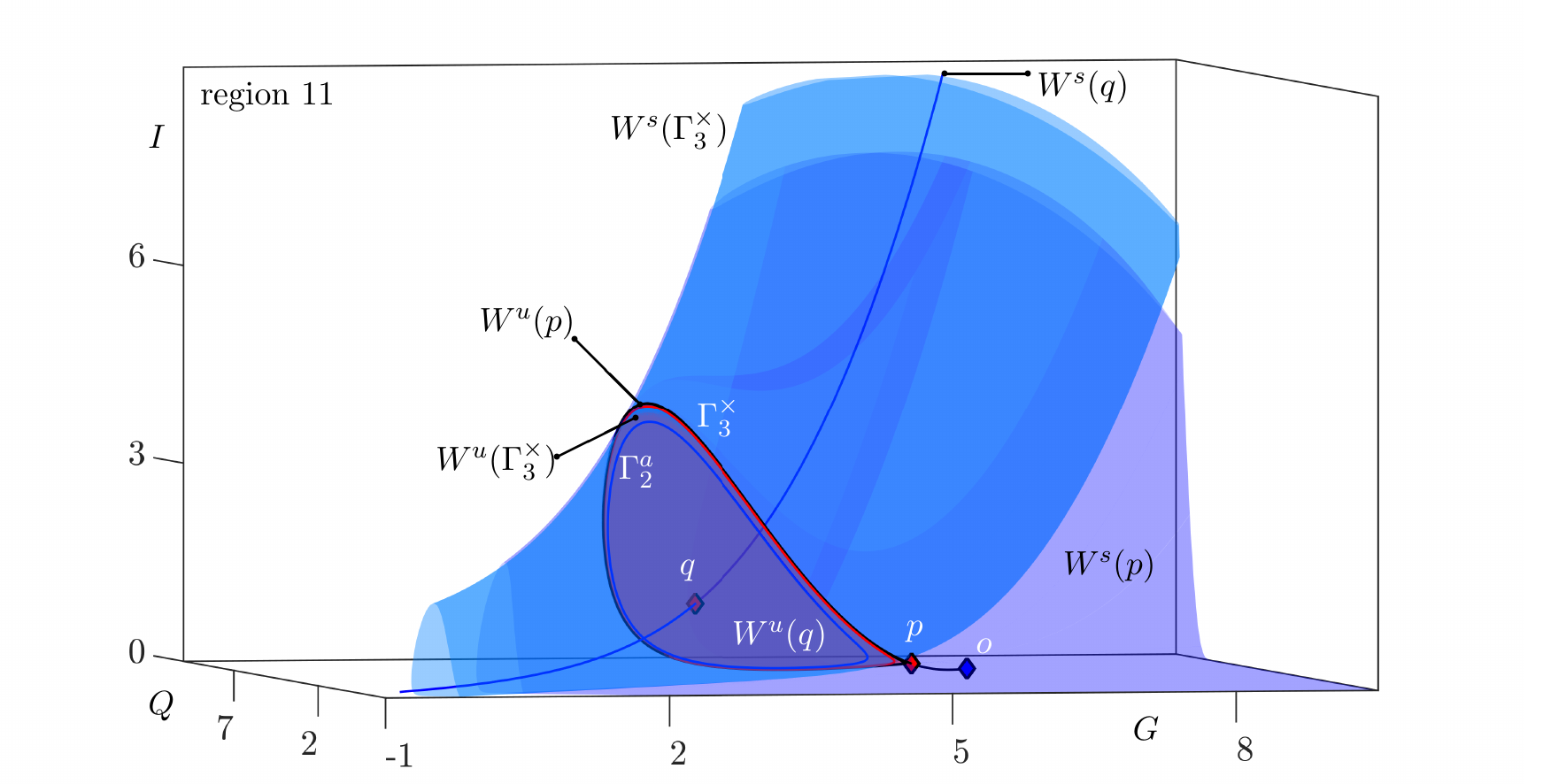}
\caption{Phase portrait in region 11 in the $(G,Q,I)$-space  showing the attracting equilibria $o$ (blue diamond), saddle equilibria $p,q$ (red diamond), attracting periodic orbit $\Gamma _{2}^{a}$ (blue curve), saddle periodic orbit $\Gamma _{3}^{\times}$ (red curve) and manifolds $W^{s}(q)$ (blue curve) and $W^{u}(p)$ (black curve). Also shown are the surfaces $W^{s}(p)$ (dark blue surface), $W^{s}(\Gamma _{3}^{\times})$ (light blue surface), $W^{u}(\Gamma _{3}^{\times})$ (red surface); see Table~\ref{regions1to11} for parameter values.}  
\label{3D_region11}
\end{figure}

\subsection{Phase portraits 10 and 11}
\label{sec:phase10to11}

Computed images of the phase portraits in new regions 10 and 11 are presented as two-dimensional projections in Fig.~\ref{stand_alone} and in the full three-dimensional phase space in Figs.~\ref{3D_region10} and ~\ref{3D_region11}.

As the sketch at the top of Fig.~\ref{Bifurcation_Sketches} and Fig.~\ref{stand_alone} show, phase portrait 10 features three attractors: the off-state  $o$, the equilibrium $q$ and the periodic orbit $\Gamma^{a}_{2}$. To fully understand the organization of phase space one must consider the phase portraits in $\mathbb{R}^3$ shown in Figs.~\ref{3D_region10}. The manifolds $W^{s}(p)$ and $W^{s}(\Gamma _{1}^{\times})$  form two different basin boundaries that delineate the three basins of attraction; hence, the system is tri-stable in the new region 10. Specifically, $W^{s}(\Gamma _{1}^{\times})$ forms a cylinder, and any initial condition on the inside of it goes to $q$; any initial condition on the outside of $W^{s}(\Gamma _{1}^{\times})$ but on the inner side of $W^{s}(p)$ goes to $\Gamma^{a}_{2}$; Finally, an initial condition on the outer side of $W^{s}(p)$ goes to $o$. Figure~\ref{3D_region10}(a2) also shows the stable manifold $W^{s}(p)$, which rolls up as a carpet-like structure accumulating on  $W^{s}(\Gamma^{\times}_{1})$.

Region 11 from region 10 differs in that, following the Hopf bifurcation H, the equilibrium $q$ is now a saddle point while the saddle periodic orbit $\Gamma^{\times}_{1}$ has disappeared; hence there are now only the two attractors $o$ and $\Gamma^{a}_{2}$ and the laser is bistable. This is already clear from the sketch at the top of Fig.~\ref{Bifurcation_Sketches} and from its computed version in Fig.~\ref{stand_alone}. The phase portraits in $\mathbb{R}^3$ of Figs.~\ref{3D_region10} show that the stable manifold $W^{s}(p)$ forms a basin boundary that delineates two basins of attraction. Notice that the surface $W^{s}(p)$ wraps around the attracting equilibrium $q$ and the periodic orbits $\Gamma^{a}_{2}$ and $\Gamma^{\times}_{3}$; as such, any initial condition on the inner side of $W^{s}(p)$ goes to the attracting periodic orbit $\Gamma^{a}_{2}$. On the other hand, an initial condition on the other side of $W^{s}(p)$ goes to the off-state $o$.

\begin{figure}[t!]
\centering
\includegraphics[width=0.93\columnwidth]{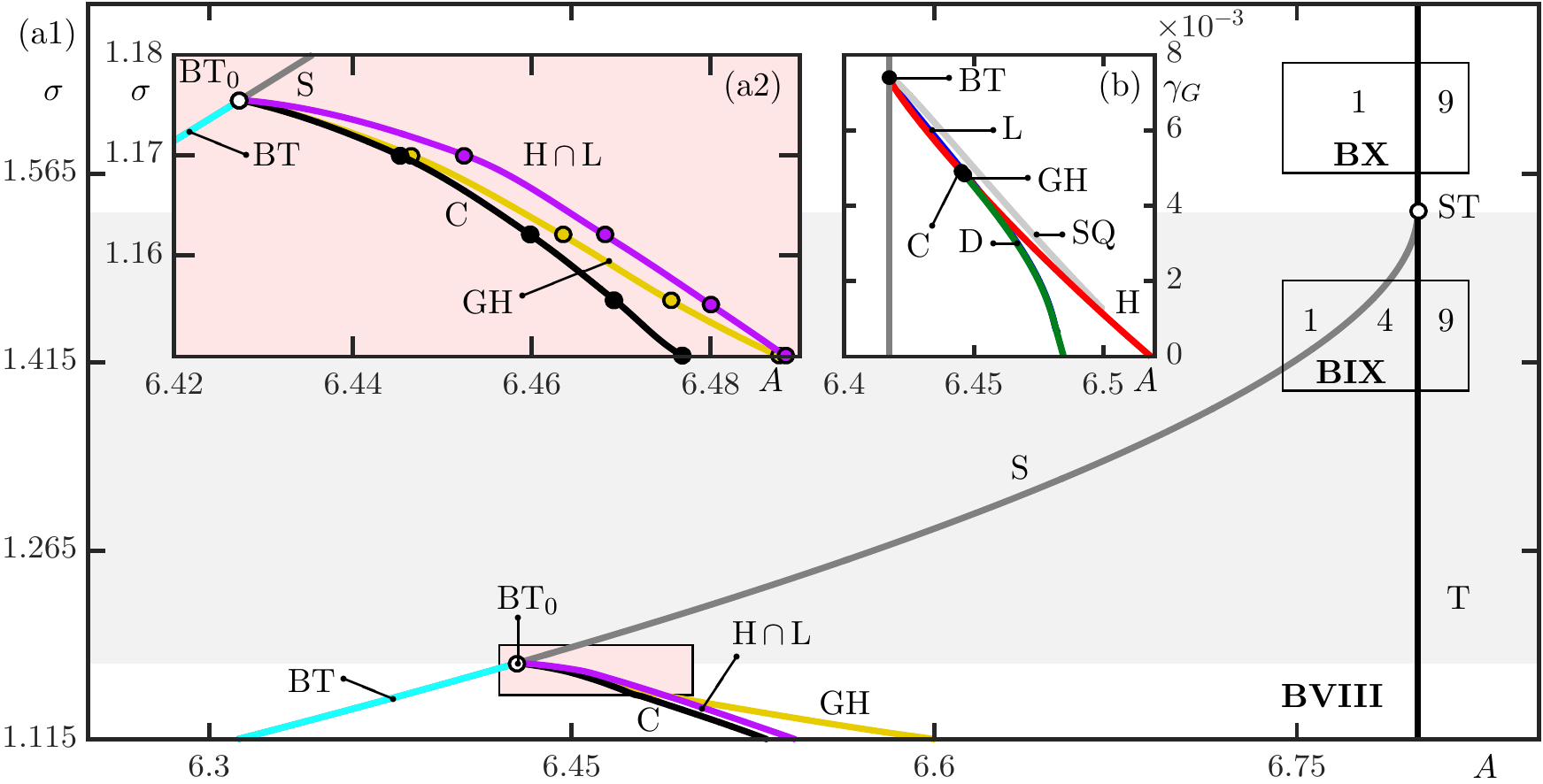}
\caption{Computed bifurcation diagram in the $(A,\sigma)$-plane. Shown in panel (a1) are the surfaces S of saddle-node (dark grey) and T of transcritical (black) bifurcations, which meet at the point ST and are curves in projection; also shown are the curves BT of Bogdanov-Takens (teal), GH of generalised Hopf (yellow), and C of cusp (black) bifurcations, and the intersection curve H$\,\cap\,$L (purple), which end at the point BT$_{0}$. The $\sigma$-range of bifurcation case BIX in between the point BT$_{0}$ and ST is shaded, and insert (a2) shows an enlargement of the shaded region near the point BT$_{0}$. Panel (b) shows the computed bifurcation diagram in the $(A,\gamma_{G})$-plane of case BVIII at $\sigma=1.1704$, with the curves S (dark grey), T (black), D (green), H (red), L (blue) and SQ (light grey) and the codimension-two points GH, C and BT; compare with Fig.~\ref{Bif_diagram_sigma}(i).}
\label{projection}  
\end{figure}

\subsection{Cases BIX and BX}
\label{sec:IXtoX}

As $\sigma$ is increased further, the point BT slides down the curve S; moreover, the bifurcation curves H, L and D move down as well toward the point where S meets the axis $\gamma_{G}=0$. This is illustrated in Fig.~\ref{projection}, where panel (a1) shows a projection onto the $(A,\sigma)$-plane; here, the codimension-one bifurcations S and T appear as curves since these bifurcations do not depend on $\gamma_{G}$. Also plotted are the projections of the codimension-two curves BT of Bogdanov-Takens bifurcation (teal), GH of generalised Hopf bifurcation (yellow) and C of cusp (black). Moreover, we plot the intersection curve H$\,\cap\,$L (purple), which is the codimension-one-plus-one event discussed earlier. The curve BT lies on S and all curves terminate at the point BT$_{0}$ where BT reaches $\gamma_{G}=0$. We determined that this happens at $\sigma^{{\rm BT}_{0}}=1.17550$ (then $\gamma_{G}=3.56084\times10^{-7}$, which is the smallest value for which the curve BT is found). Panel (a2) is an enlargement around the point BT$_{0}$; it shows how the curves BT, GH, C and H$\,\cap\,$L indeed approach BT$_{0}$ while maintaining their relative positions. The curves H, L and D have been found for $\sigma \in \left\{1.1555, 1.1623, 1.1704\right\}$, which is indicated by the dots in Fig.~\ref{projection}(a2) from which the shown curves have been rendered. Figure.~\ref{projection}(b) shows the final slice of the bifurcation diagram in the $(A,\gamma_{G})$-plane for $\sigma=1.1704$, which is at the limit of what can be achieved with the continuation software. Indeed, the bifurcation curves H, L and D are topologically as in Fig.~\ref{Bif_diagram_sigma}(i) of case BVIII; hence this is numerical confirmation that H, L and D shrink down to the point BT$_{0}$ while maintaining their relative positions, meaning that case BVIII exists all the way up to $\sigma^{{\rm BT}_{0}}=1.17550$ where the curves BT, GH, C and  H$\,\cap\,$L disappear at BT$_{0}$. For $\sigma^{{\rm BT}_{0}} < \sigma$ the new case BIX occurs, which consists of only crossing S and T for any $\gamma_{G}$, giving the transition from region 1 to region 4 to region 9. This is illustrated in Fig.~\ref{projection}(a1) by the box in the grey region indicating the $\sigma$-range of case BIX. 

For even larger $\sigma$, the curve S becomes tangent to T at the point ST at $\sigma^{\rm ST}=1.53529$, above which the curve S continues on but corresponds to a saddle node bifurcation for negative intensity $I$. Hence, for $\sigma^{\rm ST} < \sigma$ the bifurcation S is no longer physically relevant, and so not part of the bifurcation diagram, and we find the new bifurcation case BX, which involves crossing only T. Case BX involves only the transition from region 1 to region 9, and this is again illustrated in the box in the white region. This actually corresponds to the laser switching on at the laser threshold: the laser is off in region 1 and turns on in region 9 with gradually increasing intensity. This is the case for any sufficiently large value of $\sigma$ and, hence, overall, there are indeed ten bifurcation cases BI to BX.

\subsection{Expressions for the points HT$_0$, BT$_0$ and ST}
\label{BT0andST}

As we show now, the $\sigma$-values of the transition points HT$_0$, BT$_{0}$ and ST can actually be computed analytically, and their analytical values indeed agree (up to five decimal places) with the numerical ones found by continuation. Their locations can be derived from analytic expressions for equilibria of Eqs.~\eqref{Yamadamodel2} and their bifurcations; these can be obtained by extending the calculations in~\cite{Dubbeldam_Krauskopf_Self_pulsations_lasers_saturable_absorber} for $\sigma = 1$ to the case of general $\sigma \neq 1$. This was first done in~\cite{unpublishedkeyB} for a variant of the Yamada model with a different $\sigma$-dependence. 

We report here the expressions of the bifurcations of Eqs.~\eqref{Yamadamodel2} that are required to determine the points HT$_0$, BT$_0$ and ST; our calculations were performed symbolically with the help of the package Mathematica.

\medskip
\noindent
{\bf Proposition}\\
The bifurcation set of Eqs.~\eqref{Yamadamodel2} in $(A,\gamma_{G},\sigma)$-space features the following bifurcations.
\begin{itemize}
\item[(i)] 
A transcritical bifurcation T of the off-state equilibrium $o = (A,B,0)$ occurs at 
\begin{eqnarray}
A^{\rm T} = B+1 \ {\rm for \ any} \ \sigma \ {\rm and \ any} \ \gamma_G.
\label{eq:tformula}
\end{eqnarray}
For $A < A^{\rm T}$ the equilibrium $o = (A,B,0)$ is stable.\\[-2mm]
\item[(ii)] 
A saddle-node bifurcation S occurs at 
\begin{eqnarray}
A^{\rm S}(\sigma)  =  \frac{a + (B-1)\sigma + 2\sqrt{B\sigma(a - \sigma)}}{a} \ {\rm for \ any} \ \gamma_G,
\label{eq:sformula}
\end{eqnarray}
and for $ 0 < \sigma < a$.
For $A < A^{\rm S}(\sigma)$ there is only the equilibrium $o = (A,B,0)$, while for $ A^{\rm S}(B) < A$ there are also the equilibria 
\begin{eqnarray}
x_\pm & =& \left(A/(1+I_{\pm}),B\sigma/(\sigma+a I_{\pm}),I_{\pm}\right) \ {\rm with} \label{eq:equilextra}\\
I_{\pm} &=& \frac{a(A-1)-(B+1)\sigma \pm 
\sqrt{a^2 (A-1)^2 - 2 a (B A - A + B + 1) \sigma + (B+1)^2  \sigma^2}}{2a}\,.
\label{eq:iextra}
\end{eqnarray}
Note that $p = x_-$ is a saddle with two unstable eigenvalues, while $q = x_+$ may be a sink or a saddle with one unstable eigenvalue.
\item[(ii)] 
A simultaneous saddle-node and transcritical bifurcation ST occurs at 
\begin{eqnarray}
A^{\rm ST}(\sigma)=B+1 \ {\rm and} \ \sigma^{\rm ST} = \frac{a B}{B + 1} \ {\rm for \ any} \ \gamma_G.
\label{eq:stformula}
\end{eqnarray}
For $a = 1.8$ and $B=5.8$ the point ST is at $A^{\rm ST} = 6.8$ and $\sigma^{\rm ST} = 1.53529$.\\[-2mm]
\item[(iv)] 
A Hopf bifurcation H of the equilibrium $q = x_+$ occurs at
\begin{eqnarray}
\gamma^{\rm H}(\sigma)  = \frac{\sigma \left(a B - A \sigma^3\right) I_+}
     {(1 + I_+) (a I_++ \sigma) (a I_+ + \sigma (1 + \sigma + \sigma I_+))}
\label{eq:hformula}
\end{eqnarray}
where $I_+$is given by \eqref{eq:iextra}.
\item[(v)] 
A simultaneous Hopf and transcritical bifurcation HT$_0$ with $\gamma = 0$ occurs at 
\begin{eqnarray}
A^{{\rm HT}_0}(\sigma)=B+1 \ {\rm and} \ \sigma^{\rm {\rm HT}_0} = \sqrt[3]{\frac{a B}{B + 1}}.
\label{eq:ht0formula}
\end{eqnarray}
For $a = 1.8$ and $B=5.8$ the point HT$_0$ is at $A^{{\rm HT}_0} = 6.8$ and $\sigma^{{\rm HT}_0} = 1.15362$.\\[-2mm]
\item[(vi)] 
A Bogdanov-Takens bifurcation occurs at 
\begin{eqnarray}
A^{\rm BT}(\sigma)  &=&  A^{\rm S}(\sigma)  \ {\rm and} \nonumber \\
\gamma_G^{\rm BT}(\sigma) &=& 
\frac{\left(-a^2 B-a \sigma^2+ (B + 1)\sigma^3\right) \sqrt{B \sigma (a-\sigma)}                                             
+ B (a - \sigma) \sigma (a (B + 1 + \sigma - (B + 1) \sigma^2))}
{B ((B +1) \sigma - a) (a - \sigma)^2}.
\label{eq:btformula}
\end{eqnarray}
\item[(vii)] 
A Bogdanov-Takens bifurcation BT$_0$ with $\gamma_G = 0$ occurs at the point $(A^{{\rm BT}_0},\sigma^{{\rm BT}_0})$, where $A^{{\rm BT}_0} = A^{\rm S}(\sigma^{{\rm BT}_0})$ and $\sigma^{{\rm BT}_0}$ is the smallest positive real zero of 
\begin{eqnarray}
 Q(\sigma) = a^2 B - 2 a B \sigma^2 - a \sigma^3+ (B + 1) \sigma^4,
\label{eq:qpoly}
\end{eqnarray}
which exists provided that 
\begin{eqnarray}
 a > \frac{16}{27} \left(-9 B - 8 B^2 +\sqrt{27 B +108 B^2 + 144 B^3 + 64 B^4} \right).
\label{eq:qpolycond}
\end{eqnarray}
For $a = 1.8$ and $B=5.8$ the point BT$_0$ is at $A^{{\rm BT}_0} = 6.4274$ and $\sigma^{{\rm BT}_0} = 1.17550$.
\end{itemize}

\medskip
\noindent
{\bf Proof:}
The Jacobian of Eqs.~\eqref{Yamadamodel2} is 
  \begin{equation}
    \mathrm{D}F = \begin{pmatrix}
      -\gamma_{g}(I + 1)   	&  0                                    			& -G\gamma_{g}   \\
      0                             	  	&  \ \ -(\gamma_{g}(\sigma + Ia))/\sigma^2 & \ \ -(Qa\gamma_{g})/\sigma^2 \\
      I                             		&  - I                                 	 		& G - Q - 1\\
    \end{pmatrix}.
\label{eq:jacobian}
  \end{equation}
Evaluating $\mathrm{D}F$ at $o = (A,B,0)$ yields $A - B - 1$ as the eigenvalue of the eigenvector $(0,0,1)$; the other two eigenvalues are nonzero and (i) follows. The additional two equilibria \eqref{eq:equilextra} are readily found from Eqs.~\eqref{Yamadamodel2}, and they exist only when the expression under the square root in \eqref{eq:iextra} is positive. When this quadratic expression is zero then there is a saddle-node bifurcation at $x_{\rm S}$, defined as in \eqref{eq:equilextra} with $I_\pm = I_{\rm S} = \pm \sqrt{B \sigma (a -\sigma)}/a$. The formula for $A^{\rm S}(\sigma)$ in \eqref{eq:sformula} is the branch where $I_{\rm S}$ is positive and, hence, (ii) follows. Equating $A^{\rm T} = A^{\rm S}(\sigma)$ of \eqref{eq:tformula} and \eqref{eq:sformula} gives \eqref{eq:stformula} and (iii).

The condition for the Hopf bifurcation is that the Jacobian $\mathrm{D}F$ in~\eqref{eq:jacobian} evaluated at the equilibrium $q = x_+$ has purely complex conjugate eigenvalues with non-zero imaginary part. This can be verified from the characteristic polynomial, which is of the form $\chi(\lambda) = \lambda^3 -  t \lambda^2 + r \lambda - d$ (where $t$ is the trace and $d$ the determinant), by imposing the Hurwitz condition $s r - d = 0$. In this way, the expression for the Hopf bifurcation in (iv) is found. Since $I_+ \neq 0$ for $0 < \sigma$, the nominator of $\gamma^{\rm H}(\sigma)$ in \eqref{eq:hformula} has a positive zero exactly when $a B-A \sigma^3 = 0$ and \eqref{eq:ht0formula} of (v) follows with $A = A^{\rm T} = B+1$. 

The Jacobian $\mathrm{D}F$ evaluated at the saddle node $x_{\rm S}$ has a zero eigenvalue, meaning that the characteristic polynomial $\chi(\lambda)$ has zero as a root and, hence, a factor $\lambda$. The quadratic condition that zero is a double root of $\chi(\lambda)$ yields the formula for $\gamma_G^{\rm BT}(\sigma)$ in \eqref{eq:btformula} of (vi). The numerator of \eqref{eq:btformula} is of the form $V\sqrt{U} + W$, and quadratic expansion, that is, multiplication with $V \sqrt{U} - W$, yields the  eight-order polynomial 
\begin{eqnarray*}
 P(\sigma) = B \sigma (\sigma - a ) ((B + 1) \sigma - a) ((B + 1) \sigma - aB) Q(\sigma),
\end{eqnarray*}
where $Q(\sigma)$ is given by \eqref{eq:qpoly}.
The four explicit zeros $0 < a/(B + 1) < aB/(B + 1) < a$ of $ P(\sigma)$ lie outside of the range $1 < \sigma < \sigma^{\rm ST} = aB/(B + 1)$ of where the point BT$_0$ must lie. Hence, its $\sigma$-value is a zero of the forth-order polynomial $Q(\sigma)$. Its derivative is 
\begin{eqnarray*}
Q' (\sigma)= \sigma (-4 a B  - 3 a \sigma + 4 (B + 1) \sigma^2).
\end{eqnarray*}
Since its quartic term is positive, it follows that $Q(\sigma)$ has the maximum $Q(0) = a^2 B > 0$ at 0 and minima at 
\begin{eqnarray*}
\sigma_\pm = \frac{3 a \pm \sqrt{a} \sqrt{9 a+64 B+64 B^2}}{8 (B+1)}
\end{eqnarray*}
with $\sigma_- < 0 < \sigma_+$. Hence, positive roots of $Q(\sigma)$ exist when 
\begin{eqnarray*}
  Q (\sigma_+) =
    -\frac{\sqrt{a} (9 a+64 B (B+1)))^{\frac{3}{2}}+ 27 a^2+288 a B (B+1)-512 B (B+1)^2}{512 (B+1)^3}
\end{eqnarray*}
is negative. To determine the limiting case where $Q (\sigma_+) = 0$ we again perform a quadratic expansion of the numerator, which yields the second-order polynomial
\begin{eqnarray*}
R(a) = 1024 B (B + 1)^3 ( 256 B (B + 1) - 32 B (8 B + 9) a - 27 a^2);
\end{eqnarray*}
its only positive zero is 
\begin{eqnarray*}
\frac{16}{27} \left(-9 B -8 B^2 + \sqrt{27 B +108 B^2 + 144 B^3 + 64 B^4} \right),
\end{eqnarray*}
which shows \eqref{eq:qpolycond} and completes (vii). \hfill $\square$

\begin{figure}[t!]
\centering
\includegraphics[width=0.97\columnwidth]{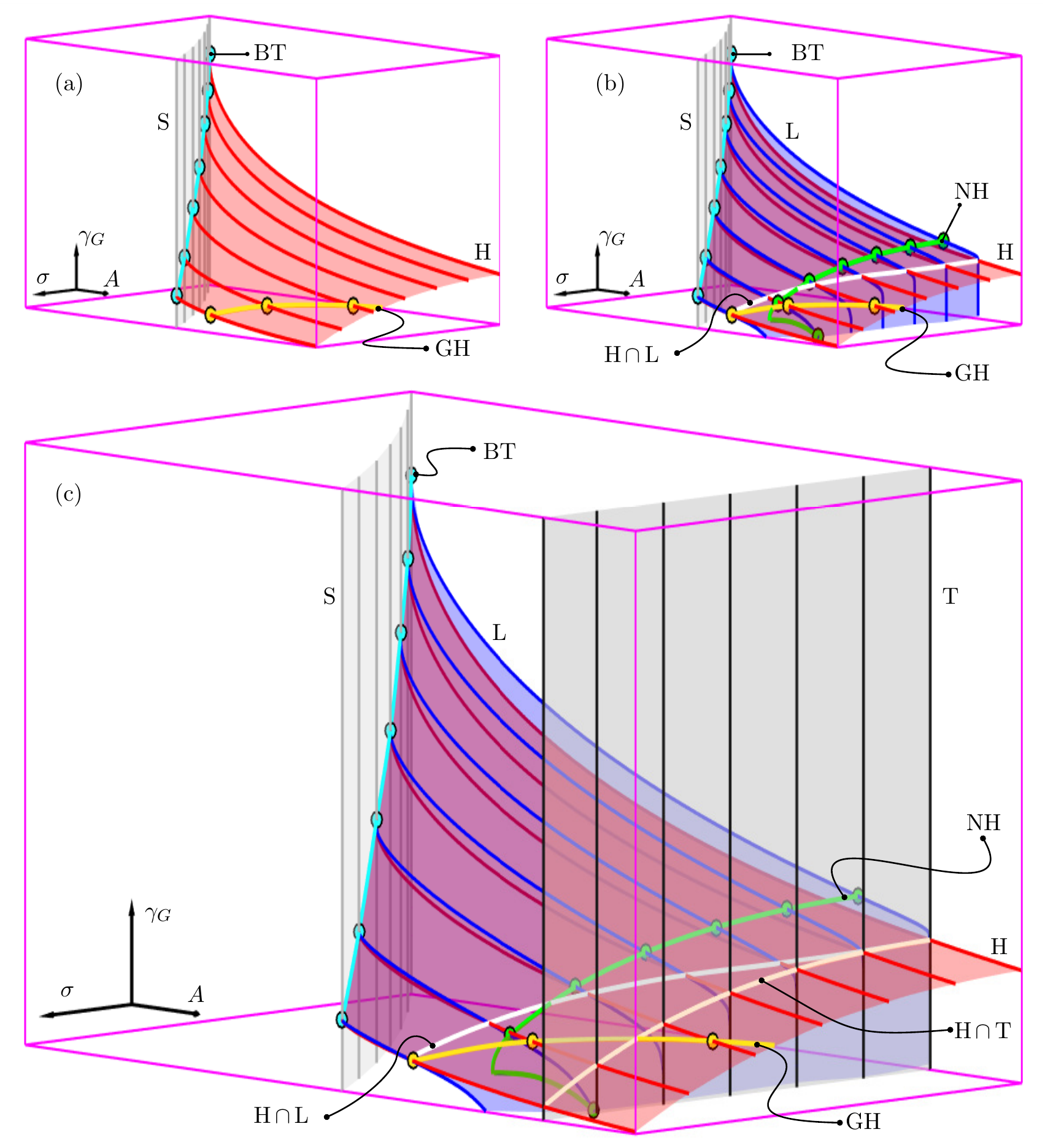}
\caption{Three parameter bifurcation diagram in the $(A,\gamma_{G},\sigma)$-space for $A \in [5.50, 7.00]$, $\gamma_{G} \in [0, 0.40]$ and $\sigma \in [0.85, 1.14]$. Panels (a) to (c) show a successive build-up of the bifurcation diagram, consisting of surfaces of codimension-one and curves of codimension-two bifurcations. The surfaces and curves are rendered from the highlighted curves and points for fixed $\sigma$.  Panel (a) shows the surfaces S of saddle-node bifurcations (light grey)  and H of Hopf bifurcations (red), which meet along the curve of Bogdanov-Takens (BT) bifurcations (teal) curve; also shown is the curve GH of generalised Hopf bifurcations (yellow). Panel (b) also shows the surface L of homoclinic bifurcations (blue), which intersects the surface H along the curve H$\,\cap\,$L (white); also shown is the curve NH of neutral saddle homoclinics (green). Panel (c) also shows the surface T of transcritical bifurcations  (dark grey) and its intersection curve H$\,\cap\,$T (cream) with H.}
\label{composite}  
\end{figure}

\subsection{Three-parameter bifurcation diagram in $(A,\gamma_{G},\sigma)$-space}
\label{sec:3Dbif}

The computed two-dimensional bifurcation diagrams in the $(A,\gamma_{G})$-plane of Figs.~\ref{Bif_diagram_sigma} and~\ref{projection} illustrate qualitative changes in shape and size of curves and regions as $\sigma$ is increased from $\sigma=1$. Figure~\ref{composite} summarises the entire transition for increasing $\sigma$ by showing the corresponding three-parameter bifurcation diagram, where $\sigma$ is the third direction. The emphasis here is not on the details of the individual cases BI--BX discussed earlier, but on providing an overview of how the loci of the different bifurcations change with $\sigma$.

Figure~\ref{composite} progressively builds up the three-parameter bifurcation diagram in $(A,\gamma_{G},\sigma)$-space from panel to panel for $A$ in the range [5.50, 7.00], $\gamma_{G}$ in the range [0, 0.40] and $\sigma$ in the range [0.85, 1.14], which spans the bifurcation cases BI to BVI. Bifurcation cases beyond case BVI are very small and would be hard to decern on the scale shown. In this representation, codimension-one bifurcations are now surfaces while codimension-two bifurcations are now curves. Panel (a) shows the surface S of saddle-node bifurcations (light grey) which has the shape of a curved vertical sheet. It meets the surface H of Hopf bifurcations (red) along the curve BT of  Bogdanov-Takens bifurcations (teal); moreover, on H there is the curve GH of general Hopf  bifurcations (yellow). Notice how the surface H dips downward towards $\gamma_{G}=0$ as $\sigma$ increases while the curves BT and GH on H come down towards smaller $\gamma_{G}$. Figure~\ref{composite}(b) shows also the surface L of homoclinic bifurcations (blue) that lies above the surface H to the immediate right of S from where H and L emanate. However, L goes below the surface H after intersecting it along the curve H$\,\cap\,$L (white). This curve is not a codimension-two bifurcation curve but rather corresponds to a codimension-one-plus-one event where two codimension-one bifurcations happen to occur simultaneously. Also shown is the curve NH of neutral-saddle homoclinics (green) on L; it crosses H, has a fold point with respect to $\sigma$ (at the point FNH) and then ends at $\gamma_{G}=0$. The final panel (c) now shows the entire bifurcation diagram in the $(A,\gamma_{G},\sigma)$-space. Also shown is the surface T of transcritical bifurcations (dark grey) which is a vertical plane to the right of the surface L; notice how T intersects H along the curve H$\,\cap\,$T (cream), which is a codimension-one-plus-one event.

Figure~\ref{composite} is designed to give detailed insight into the relative positions of surfaces; therefore, these have been rendered in a transparent fashion so that changes in the general positions of curves and surfaces are easily seen. Although this compact representation does not show all the details that have been discussed at great length in the previous section, it provides an overall  impression of the regions and surfaces in the $(A,\gamma_{G},\sigma)$-space. In particular, all the surfaces come down towards smaller values of $\gamma_{G}$ as $\sigma$ is increased. This shows that the regions bounded by these surfaces also become smaller and can only be found close to $\gamma_{G}=0$ for larger $\sigma$.

\newpage 

\section{Conclusions}

We presented a classification of the dynamics of a self-pulsing laser with saturable absorber as described by the Yamada model. More specifically, for the most relevant case with complicated dynamics organised by a Bogdanov-Takens point, we presented the bifurcation structure that explains how the observable dynamics depends on the pump parameter A, the gain decay parameter $\gamma_{G}$, and the decay ratio $\sigma$ between $\gamma_{G}$ the decay parameter $\gamma_{Q}$ in the absorber section. The new aspect of this work is the additional dependence on the parameter $\sigma$. Our results expand the previous work in~\cite{Dubbeldam_Krauskopf_Self_pulsations_lasers_saturable_absorber} that only considered the case of equal decay times. As a first result, we showed here that the respective bifurcation diagram, BI in our notation, is actually valid for all $0<\sigma \leq 1$, that is, for any laser with saturable absorber with $\gamma_{G} \leq \gamma_{Q}$; this explains why the initial study in~\cite{Dubbeldam_Krauskopf_Self_pulsations_lasers_saturable_absorber} has been more widely applicable than initially thought. Morever, our study gives a complete understanding of the dynamics for any value of $\sigma$. In fact, for $1 < \sigma$ such that $\gamma_{Q} < \gamma_{G}$, we found new dynamics and changes in the composition of the bifurcation diagram, including two new regions appearing and other regions disappearing. We provided sketches of generic bifurcation cases BI to BX and transitions between them via codimension-three events as $\sigma$ increases. These sketches are accompanied by evidence in the $(A, \gamma_{G})$-plane of careful numerical investigations, including enlargements of regions that verify the presence or absence of certain dynamics. While our bifurcation study is complete from a mathematical perspective, some regions are very small and can be found only by very careful numerical investigation. In fact, we do not expect all regions to be physically accessible in an experiment. On the other hand, the overall dynamics presented is of immediate relevance, and it is our hope that our results may serve as a ``road map" for future experimental investigations.

One concrete example of an observable prediction from our work is the following. For bifurcation cases BV to BX, for small$\gamma_G$ the homoclinic bifurcation L is no longer immediately followed by the transcritical bifurcation T. In practical terms this means that self-pulsations will be observable for suitable initial conditions prior to the threshold of the laser; this is in contrast to cases BI to BIV, where the onset of self-pulsations effectively coincides with the laser threshold~\cite{Dubbeldam_Krauskopf_Self_pulsations_lasers_saturable_absorber}. The transition between these two scenarios is due to a neutral saddle homoclinic bifurcation point existing for $\gamma_G = 0$, which is a special point of codimension three; in spite of the fact that it represents a global bifurcation, it might be possible to study this point by blow-up techniques similar to those in~\cite{huber}.

In addition to the bifurcation analysis, we have also looked in detail at the geometry of invariant objects in the full three-dimensional phase space. For the different regions of qualitative dynamics, we presented complete computed images of the respective phase portraits. In particular, we determined and discussed in detail the geometry of the two-dimensional stable manifolds of the relevant invariant objects, equilibria and periodic orbits of saddle type, that form the boundaries between different basins of attraction. In this way, we clarified the organisation of multi-stabilty, as well as the exact nature of the excitability threshold in this laser system. These results are made possible by employing advanced numerical techniques, specifically the continuation of solutions of suitable two-point boundary value problems~\cite{krauskopf2007numerical}. Showing global manifolds as surfaces in the three-dimensional phase space is important from a theoretical as well as a practical perspective to determine the nature and number of attractors, their locations in phase space and the locations of the boundaries of their basins. Although, we give renditions of all phase portraits at specific values of the system parameters, they nevertheless provide a general representation of what the organisation of phase space looks like in the particular region. Our representative phase portraits have the potential to predict effects that may appear to be counter-intuitive. An example is the observation that an excitable response may not occur if the system is perturbed with an input pulse that is too large. This is explained by the geometry of the excitability threshold, formed by the stable manifold of a saddle equilibrium in the relevant region, which folds over in the region of the off-state to create an upper bound of intensity perturbations that result in an output pulse.

We expect that our work will be of direct interest to any study of laser devices with different decay times of gain and absorber, which is the typical situation. An example is the characterisation of a micropillar laser with an intra-cavity saturable absorber, which was shown in~\cite{barbay2011excitability, selmi2016spike} to be described well by the Yamda model with $\sigma$ close to one. Subsequent work in~\cite{terrien2017asymmetric, terrien2018pulse} considered sustaining  pulse trains in such a device via optical feedback, where the gain decays faster than the absorption. In recent work~\cite{otupiri2018experimental}, we consider a relatively short all-fibre laser, comprised of a 1~m gain section and an absorber section ranging from 0.21~m to 1.48~m, and showed that it displays Q-switching. Moreover, for this laser system pulsing mechanism and threshold shows good agreement with the predictions of the Yamada model. In particular, the onset of oscillations and subsequent pulse characteristics, either via a homoclinic bifurcation or  a Hopf bifurcation, agrees very well with experimental observations. For comparison with a specific experiment, it is necessary to translate the physical parameters of the device to associated parameters of the model; these include the input pump power to $A$ and decay life times $\gamma_{G}$ and $\gamma_{Q}$ of gain and absorber medium, respectively. This can be a challenge in light of experimental difficulties to determine the exact values of physical parameters. Nevertheless, the work presented here can serve as a guide, in particular, since certain dynamics of interest, such as multi-stabilty and excitability, exists over reasonably large ranges of the system parameters.

\section*{Ackknowlegdements}
We thank Andrus Giraldo, Soizic Terrien and Stefan Ruschel for helpful discussions.



\bibliographystyle{unsrt}
\bibliography{okb_decay_ratio_ref}

\begin{thebibliography}{10}

\bibitem{abraham1988dynamical}
N.~B. Abraham, P.~Mandel, and L.~M. Narducci.
\newblock I dynamical instabilities and pulsations in lasers.
\newblock In {\em Progress in Optics}, volume~25, pages 1--190. Elsevier, 1988.

\bibitem{velikhov1983pulsed}
E.~P. Velikhov, V.~Yu. Baranov, V.~S. Letokhov, E.~A. Ryabov, and A.~N.
  Starostin.
\newblock Pulsed {CO}$_2$-lasers and their application for isotope separation.
\newblock 1983.

\bibitem{heard1963ultra}
H.~G. Heard.
\newblock Ultra-violet gas laser at room temperature.
\newblock {\em Nature}, 200(4907):667, 1963.

\bibitem{little1999metal}
C.~E. Little.
\newblock Metal vapour lasers: Physics, engineering and applications.
\newblock {\em Metal Vapour Lasers: Physics, Engineering and Applications, by
  Christopher E. Little, pp. 646. ISBN 0-471-97387-4. Wiley-VCH, March 1999.},
  page 646, 1999.

\bibitem{magni1986resonators}
V.~Magni.
\newblock Resonators for solid-state lasers with large-volume fundamental mode
  and high alignment stability.
\newblock {\em Applied Optics}, 25(1):107--117, 1986.

\bibitem{yu2013topological}
H.~Yu, H.~Zhang, Y.~Wang, C.~Zhao, B.~Wang, S.~Wen, H.~Zhang, and J.~Wang.
\newblock Topological insulator as an optical modulator for pulsed solid-state
  lasers.
\newblock {\em Laser \& Photonics Reviews}, 7(6):L77--L83, 2013.

\bibitem{keller1996semiconductor}
U.~Keller, K.~J. Weingarten, F.~X. Kartner, D.l Kopf, B.~Braun, I.~D. Jung,
  R.~Fluck, C.~Honninger, N.~Matuschek, and J.~A. Der~Au.
\newblock Semiconductor saturable absorber mirrors (sesam's) for femtosecond to
  nanosecond pulse generation in solid-state lasers.
\newblock {\em IEEE Journal of Selected Topics in Quantum Electronics},
  2(3):435--453, 1996.

\bibitem{rulliere2005femtosecond}
C.~(Ed.) Rulliere.
\newblock {\em Femtosecond laser pulses}.
\newblock Springer, 2005.

\bibitem{tam1986applications}
A.~C. Tam.
\newblock Applications of photoacoustic sensing techniques.
\newblock {\em Reviews of Modern Physics}, 58(2):381, 1986.

\bibitem{novak1994millimetre}
D.~Novak and R.~S. Tucker.
\newblock Millimetre-wave signal generation using pulsed semiconductor lasers.
\newblock {\em Electronics Letters}, 30(17):1430--1431, 1994.

\bibitem{duarte1990dye}
F.~J. Duarte, P.~Kelley, L.~W. Hillman, and P.~F. Liao.
\newblock {\em Dye laser principles: with applications}.
\newblock Academic Press, 1990.

\bibitem{hamlin1991high}
S.~J. Hamlin, J.~D. Myers, and M.~J. Myers.
\newblock High repetition rate q-switched erbium glass lasers.
\newblock In {\em Eyesafe Lasers: Components, Systems, and Applications},
  volume 1419, pages 100--106. International Society for Optics and Photonics,
  1991.

\bibitem{gibson1993lasers}
K.~F. Gibson and W.~G. Kernohant.
\newblock Lasers in medicine-a review.
\newblock {\em Journal of Medical Engineering \& Technology}, 17(2):51--57,
  1993.

\bibitem{svelto1998principles}
O.~Svelto and D.~C. Hanna.
\newblock {\em Principles of lasers}, volume~4.
\newblock Springer, 1998.

\bibitem{siegman1971introduction}
A.~E. Siegman and A.~E. Siegman.
\newblock {\em An introduction to lasers and masers}, volume 122.
\newblock McGraw-Hill, New York, 1971.

\bibitem{erneux1988q}
T.~Erneux.
\newblock Q-switching bifurcation in a laser with a saturable absorber.
\newblock {\em JOSA B}, 5(5):1063--1069, 1988.

\bibitem{haus2000mode}
H.~A. Haus.
\newblock Mode-locking of lasers.
\newblock {\em IEEE Journal of Selected Topics in Quantum Electronics},
  6(6):1173--1185, 2000.

\bibitem{smith1970mode}
P.~W. Smith.
\newblock Mode-locking of lasers.
\newblock {\em Proceedings of IEEE}, 58(9):1342--1357, 1970.

\bibitem{yamada_theoretical_1993}
M.~Yamada.
\newblock A theoretical analysis of self-sustained pulsation phenomena in
  narrow-stripe semiconductor lasers.
\newblock {\em IEEE Journal of Quantum Electronics}, 29(5):1330--1336, 1993.

\bibitem{ueno1985conditions}
M.~Ueno and R.~Lang.
\newblock Conditions for self-sustained pulsation and bistability in
  semiconductor lasers.
\newblock {\em Journal of Applied Physics}, 58(4):1689--1692, 1985.

\bibitem{Dubbeldam_Krauskopf_Self_pulsations_lasers_saturable_absorber}
J.~L.~A. Dubbeldam and B.~Krauskopf.
\newblock Self-pulsations of lasers with saturable absorber: dynamics and
  bifurcations.
\newblock {\em Optics Communications}, 159(4):325--338, 1999.

\bibitem{peterson1999dynamics}
P.~Peterson, A.~Gavrielides, M.~P. Sharma, and T.~Emeux.
\newblock Dynamics of passively q-switched microchip lasers.
\newblock {\em IEEE Journal of Quantum Electronics}, 35(8):1247--1256, 1999.

\bibitem{erneux2000pulse}
T.~Erneux, P.~Peterson, and journal={The European Physical Journal D: Atomic,
  Molecular, Optical and Plasma Physics} volume={10} number={3}
  pages={423--431} year={2000}~publisher={Springer} Gavrielides, A.
\newblock The pulse shape of a passively q-switched microchip laser.

\bibitem{barbay2011excitability}
S.~Barbay, R.~Kuszelewicz, and A.~M. Yacomotti.
\newblock Excitability in a semiconductor laser with saturable absorber.
\newblock {\em Optics Letters}, 36(23):4476--4478, 2011.

\bibitem{selmi2016spike}
F.~Selmi, R.~Braive, G.~Beaudoin, I.~Sagnes, T.~Kuszelewicz, R.and~Erneux, and
  S.~Barbay.
\newblock Spike latency and response properties of an excitable micropillar
  laser.
\newblock {\em Physical Review E}, 94(4):042219, 2016.

\bibitem{shastri2013graphene}
B.~J. Shastri, M.~A. Nahmias, A.~N. Tait, Y.~Tian, B.~Wu, and P.~R. Prucnal.
\newblock Graphene excitable laser for photonic spike processing.
\newblock In {\em 2013 IEEE Photonics Conference}, pages 1--2. IEEE, 2013.

\bibitem{otupiri2018experimental}
R.~Otupiri, B.~Garbin, B.~Krauskopf, and N.~G.~R. Broderick.
\newblock Experimental and numerical characterization of an all-fiber laser
  with a saturable absorber.
\newblock {\em Optics Letters}, 43(20):4945--4948, 2018.

\bibitem{GH}
J.~Guckenheimer and P.~Holmes.
\newblock {\em Nonlinear oscillations, dynamical systems, and bifurcation of
  vector fields}.
\newblock 1983.

\bibitem{Kuz}
Yu.~A. Kuznetsov.
\newblock {\em Elements of Applied Bifurcation Theory}.
\newblock Springer-Verlag, New York, 3rd edition, 2004.

\bibitem{unpublishedkeyB}
M.~Bosewitz, B.~Krauskopf, and N.~G.~R. Broderick.
\newblock Bifurcations for two different $\gamma$.
\newblock Technical report, 2015.

\bibitem{huber}
A.~Huber and P.~Szmolyan.
\newblock Geometric singular perturbation analysis of the yamada model.
\newblock {\em SIAM J. Applied Dynamicl Systems}, 4(3):607--648, 2005.

\bibitem{krauskopf2007numerical}
B.~Krauskopf, H.~M. Osinga, and J.~(Eds.) Gal{\'a}n-Vioque.
\newblock {\em Numerical continuation methods for dynamical systems}.
\newblock Springer, 2007.

\bibitem{terrien2017asymmetric}
S.~Terrien, B.~Krauskopf, N.~G.~R. Broderick, L.~Andr{\'e}oli, F.~Selmi,
  R.~Braive, G.~Beaudoin, I.~Sagnes, and S.~Barbay.
\newblock Asymmetric noise sensitivity of pulse trains in an excitable
  microlaser with delayed optical feedback.
\newblock {\em Physical Review A}, 96(4):043863, 2017.

\bibitem{terrien2018pulse}
S.~Terrien, B.~Krauskopf, N.~G.~R. Broderick, R.~Braive, G.~Beaudoin,
  I.~Sagnes, and S.~Barbay.
\newblock Pulse train interaction and control in a microcavity laser with
  delayed optical feedback.
\newblock {\em Optics Letters}, 43(13):3013--3016, 2018.

\end{thebibliography}

\end{document}